\documentclass[11pt]{amsart}
\usepackage{amsmath}
\usepackage{a4wide}
\usepackage[utf8]{inputenc}
\usepackage{amssymb}
\usepackage{amsopn}
\usepackage{epsfig}
\usepackage{amsfonts}
\usepackage{latexsym}
\usepackage{amsthm}
\usepackage{enumerate}
\usepackage[UKenglish]{babel}
\usepackage{verbatim}
\usepackage{color}
\usepackage{calrsfs}
\usepackage{pgf}
\usepackage{tikz}
\usetikzlibrary{positioning} 
\usepackage{xcolor} 

\usepackage{mathrsfs}   
\usepackage{bm}

 \usetikzlibrary{arrows,automata}

\DeclareMathAlphabet{\pazocal}{OMS}{zplm}{m}{n}
\newtheorem{theorem}{Theorem}[section]
\newtheorem{lemma}[theorem]{Lemma}
\newtheorem{proposition}[theorem]{Proposition}

\theoremstyle{definition}
\newtheorem{definition}[theorem]{Definition}
\newtheorem{example}[theorem]{Example}

\theoremstyle{remark}
\newtheorem{remark}[theorem]{Remark}

\numberwithin{equation}{section}

\newcommand{\R}{\ensuremath{\mathbb{R}}}
\newcommand{\N}{\ensuremath{\mathbb{N}}}

\newcommand{\+}{ {\ensuremath{\oplus}}}

\renewcommand{\c}{ {\mathbf{c}}}
\renewcommand{\d}{ {\mathbf{d}}}

\renewcommand{\S}{\ensuremath{\pazocal{S}}}
\renewcommand{\t}{ {\mathbf{t}}}

\renewcommand{\u}{\ensuremath{\pazocal{U}}}

\newcommand{\us}{\mathbf{U}}

\newcommand{\U} {{\mathbf{U}}}

\newcommand{\su}{ {\mathbf{u}}}

\newcommand{\Om}{ {\Omega}}

\newcommand{\set}[1]{\left\{#1\right\}}

\newcommand{\la}{\lambda}
\newcommand{\ga}{\gamma}

\newcommand{\f}{\infty}
\newcommand{\de}{\delta}

\newcommand{\al}{\alpha}
\newcommand{\lle}{\preccurlyeq}
\newcommand{\lge}{\succcurlyeq}
\renewcommand{\a}{ \mathbf{a}}
\newcommand{\si}{\sigma}

\newcommand{\ra}{\rightarrow}

\begin{document}

\title{Periodic unique codings of  fat Sierpinski gasket}

\author{Derong Kong}
\address[D. Kong]{College of Mathematics and Statistics, Center of Mathematics, Chongqing University, Chongqing, 401331, P.R.China.}
\email{derongkong@126.com}

\author{Yuhan Zhang}
\address[Y. Zhang]{College of Mathematics and Statistics,   Chongqing University, Chongqing, 401331, P.R.China}
\email{zyh240207@163.com}

\dedicatory{}



\begin{abstract}
For $\beta>1$ let $S_\beta$ be the Sierpinski gasket generated by the iterated function system
\[\left\{f_{\alpha_0}(x,y)=\Big(\frac{x}{\beta},\frac{y}{\beta}\Big), \quad f_{\alpha_1}(x,y)=\Big(\frac{x+1}{\beta}, \frac{y}{\beta}\Big), \quad f_{\alpha_2}(x,y)=\Big(\frac{x}{\beta}, \frac{y+1}{\beta}\Big)\right\}.\]
 If  $\beta\in(1,2]$, then the overlap region $O_\beta:=\bigcup_{i\ne j}f_{\alpha_i}(\Delta_\beta)\cap f_{\alpha_j}(\Delta_\beta)$ is nonempty, where $\Delta_\beta$ is the convex hull of $S_\beta$. In this paper we study the periodic codings of the univoque set
\[
\mathbf U_\beta:=\left\{(d_i)_{i=1}^\f\in\{(0,0), (1,0), (0,1)\}^\N: \sum_{i=1}^\f d_{n+i}\beta^{-i}\in S_\beta\setminus O_\beta~\forall n\ge 0\right\}.
\]
More precisely, we determine for each $k\in\N$ the smallest base $\beta_k\in(1,2]$ such that for any $\beta>\beta_k$ the set $\mathbf U_\beta$ contains a   sequence of smallest period $k$. We show that each $\beta_k$ is a Perron number, and the sequence $(\beta_k)$ has infinitely many accumulation points. Furthermore, we show that $\beta_{3k}>\beta_{3\ell}$ if and only if $k$ is larger than $\ell$ in the Sharkovskii ordering; and the sequences $ (\beta_{3\ell+1}), (\beta_{3\ell+2})$  decreasingly converge to  the same limit point $\beta_a\approx 1.55898$, respectively.  In particular, we find that $\beta_{6m+4}=\beta_{3m+2}$ for all $m\ge 0$. Consequently, we prove that if $\mathbf U_\beta$ contains a sequence of smallest period $2$ or $4$, then $\mathbf U_\beta$ contains a sequence of smallest period $k$ for any $k\in\N$.

  \end{abstract}

\maketitle

\section{Introduction}\label{s1}

Given $\beta\in(1,2]$,  let  $S_\beta$ be the  Sierpinski gasket  in $\R^2$  generated by the \emph{iterated function system} (IFS)
 \[
 f_{d}(x)=\frac{x+d}{\beta}\quad\textrm{with}\quad d\in \Omega:=\set{\al_0,\al_1,\al_2},
 \]
 where $\al_0=(0,0), \al_1=(1,0)$ and $\al_2=(0,1)$.
 Then $S_\beta$ is the unique non-empty compact set in $\R^2$ satisfying
 $S_\beta=\bigcup_{d\in \Omega} f_{d}(S_\beta)$, which can be written as
  \begin{equation}\label{eq:S-beta}
S_\beta=\set{\sum_{i=1}^\f\frac{d_i}{\beta^i}: d_i\in \Omega\textrm{ for all }i\ge 1}.
 \end{equation}
   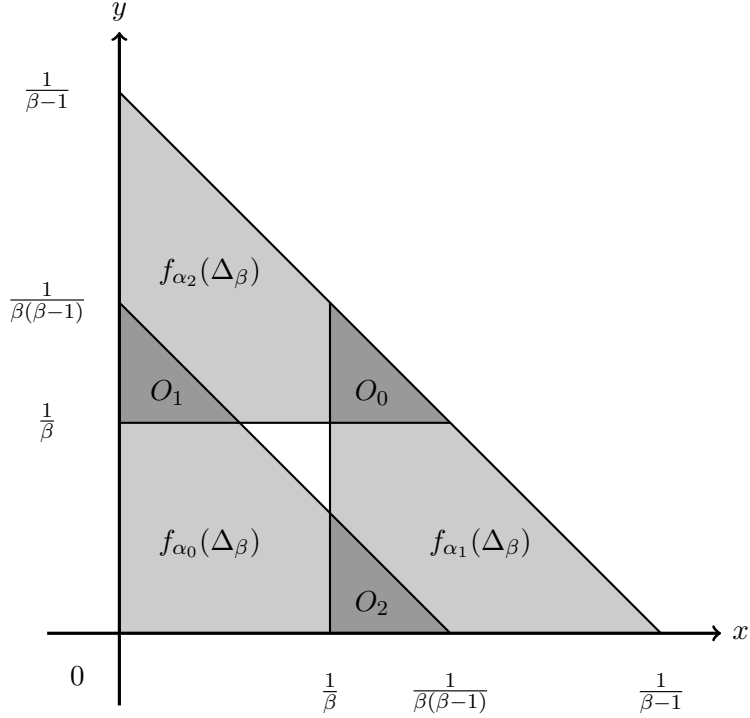
\begin{figure}[h!]
\begin{center}
\begin{tikzpicture}[
    scale=8,
    axis/.style={very thick, ->},
    important line/.style={thick},
    dashed line/.style={dashed, thin},
    pile/.style={thick, ->, >=stealth', shorten <=2pt, shorten
    >=2pt},
    every node/.style={color=black}
    ]


      \node[] at (-0.12,0.9){$\frac{1}{\beta-1}$};
         \node[] at (-0.12,0.55){$\frac{1}{\beta(\beta-1)}$};
         \node[] at (-0.12,0.35){$\frac{1}{\beta}$};

       \fill[black!40](0,0.35)--(0.2,0.35)--(0,0.55)--cycle;
         \fill[black!40](0.35,0)--(0.55,0)--(0.35,0.2)--cycle;
      \fill[black!40](0.35,0.35)--(0.55,0.35)--(0.35,0.55)--cycle;

      \fill[black!20](0,0)--(0.35,0)--(0.35,0.2)--(0.2,0.35)--(0,0.35)--cycle;

            \fill[black!20](0.55,0)--(0.9,0)--(0.55,0.35)--(0.35,0.35)--(0.35,0.2)--cycle;

                  \fill[black!20](0,0.55)--(0.2,0.35)--(0.35,0.35)--(0.35,0.55)--(0,0.9)--cycle;
             \draw[axis] (-0.12,0)  -- (1.0,0) node(xline)[right]
        {$x$};
    \draw[axis] (0,-0.12) -- (0,1.0) node(yline)[above] {$y$};

    \node[] at (-0.07,-0.07){$0$};

                    \draw[important line]  (0,0.9)--(0.9,0);
    \draw[important line] (0,0.55)--(0.55,0);
    \draw[important line] (0.35, 0)--(0.35,0.55);

       \draw[important line] (0, 0.35)--(0.55,0.35);

   \node[] at (0.9, -0.1){$\frac{1}{\beta-1}$};

     \node[] at (0.55, -0.1){$\frac{1}{\beta(\beta-1)}$};
      \node[] at (0.35, -0.1){$\frac{1}{\beta}$};

   \node[] at (0.42, 0.05){$O_2$};
      \node[] at (0.42, 0.4){$O_0$};
         \node[] at (0.08, 0.4){$O_1$};

        \node[] at (0.15, 0.15){$f_{\al_0}(\Delta_\beta)$};
          \node[] at (0.15, 0.6){$f_{\al_2}(\Delta_\beta)$};
              \node[] at (0.6, 0.15){$f_{\al_1}(\Delta_\beta)$};
\end{tikzpicture}
\end{center}
\caption{The graph of the first generation  of $S_\beta$ with $\beta=18/11\approx 1.63636$.  Then $f_{\al_0}(\Delta_\beta)$ is the left-bottom triangle, $f_{\al_1}(\Delta_\beta)$ is the right-bottom triangle, and $f_{\al_2}(\Delta_\beta)$ is the top triangle. The overlap region is  $O_\beta=O_0\cup O_1\cup O_2$.}\label{fig:1}
\end{figure}

  Clearly, the convex hull $\Delta_\beta$ of $S_\beta$ is a  right triangle with three vertices $(0,0), (\frac{1}{\beta-1},0)$ and $(0,\frac{1}{\beta})$.  Note by (\ref{eq:S-beta}) that for each  $x\in\S_\beta$ there exists a sequence $(d_i)\in \Omega^\N$ such that
  $
 x=\sum_{i=1}^\f\frac{d_i}{\beta^i},
 $
and the infinite sequence $(d_i)\in \Omega^\N$ is called a \emph{coding} of $x$ with respect to the \emph{alphabet} $\Omega$. Since  $\beta\in(1,2]$, a point in $S_\beta$ may have multiple codings. In particular, when $\beta\in(1,3/2]$ we have $S_\beta=\Delta_\beta$. In this case, Lebesgue almost every $x\in S_\beta$ has a continuum of codings (see \cite[Theorem 3.5]{Sidorov_2007}). However, when $\beta\in(3/2,2]$ the fat Sierpinski gasket $S_\beta$ has more complicated structure and has attracted much attention in the past twenty years (cf.~\cite{Jordan-05, Jordan-Pollicott-06}). {Broomhead, Montaldi and Sidorov \cite{Bro-Mon-Sid-04} proved that for $\beta\le \beta_*\approx 1.54369$ the self-similar set $S_\beta$ has nonempty interior; and for $\beta>\sqrt{3}$ the set $S_\beta$ has zero Lebesgue measure, where $\beta_*$ is a zero of $x^3-2x^2+2x-2$.} Simon and Solomyak \cite{Simon-Solomyak-02} showed that the Hausdorff dimension $\dim_H S_\beta<\log 3/\log\beta$ for a dense set of $\beta\in(3/2, 2]$, where $\log 3/\log\beta$ is the self-similarity dimension of $S_\beta$. On the other hand, Hochman \cite[Theorem 1.16]{Hochman-2015} showed that $\dim_H S_\beta=\log 3/\log\beta$ for all $\beta\in(3/2, 2]$ up to a set of zero packing dimension.

Recently,  Li and the first author   \cite{Kong-Li-2020} considered the set of points in $S_\beta$ with a unique coding. In particular, they introduced the following set
\[
U_\beta:=\set{\sum_{i=1}^\f\frac{d_i}{\beta^i}\in S_\beta: ~\sum_{i=1}^\f\frac{d_{n+i}}{\beta^i}\,\notin\, O_\beta\quad\textrm{for all }n\ge 0},
\]
where
$
O_\beta:=\bigcup_{c,d\in \Omega, c\ne d}f_{c}(\Delta_\beta)\cap f_{d}(\Delta_\beta)$ (see Figure \ref{fig:1}).
Note that in \cite{Kong-Li-2020} the set $U_\beta$ was called an \emph{intrinsic univoque set}  and  denoted by $\widetilde{\u_\beta}$. It is known that each point in $U_\beta$ has a unique coding. In \cite{Kong-Li-2020} they showed that there exists a transcendental number $\beta_c\approx 1.55263$, {which will be define in (\ref{eq:beta-c}) below,}  such that $U_\beta$ has positive Hausdorff dimension  if and only if $\beta>\beta_c$.

 We   recall  the \emph{Sharkovskii ordering} for the natural numbers  defined as follows (cf.~\cite{Sharkovski-1964}):
\begin{equation}\label{eq:sharkovski}
\begin{array}
  {cccccccccccc}
      &3&\rhd&5&\rhd&7&\rhd&\ldots&\rhd&2m+1&\rhd&\ldots\\
  \rhd&2\cdot 3&\rhd&2\cdot 5&\rhd&2\cdot 7&\rhd&\ldots&\rhd&2(2m+1)&\rhd&\ldots\\
  \rhd&2^2\cdot 3&\rhd&2^2\cdot 5&\rhd&2^2\cdot 7&\rhd&\ldots&\rhd&2^2(2m+1)&\rhd&\ldots\\
  &\vdots&&\vdots&&\vdots&&&&\vdots&&\\
  \rhd&2^n\cdot 3&\rhd&2^n\cdot 5&\rhd&2^n\cdot 7&\rhd&\ldots&\rhd&2^n(2m+1)&\rhd&\ldots\\
  &\vdots&&\vdots&&\vdots&&&&\vdots&&\\
  &&&\ldots&\rhd&8&\rhd&4&\rhd&2&\rhd&1.
\end{array}
\end{equation}
Here $k\rhd \ell$ means $k$ is {larger than  $\ell$ in the Sharkovskii ordering}. The famous Sharkovskii's theorem says that for a continuous map $f: \R\to\R$, if $k\rhd\ell$ in Sharkovskii ordering and if $f$ has a point of smallest period $k$, then $f$ has a point of smallest period $\ell$. Recently, Allouche et al.~\cite{Allouche-Clarke-Sidorov-2009} found the intimate connection between   Sharkovskii ordering and periodic unique beta expansions.

Since each point in $U_\beta$ has a unique coding, there is a bijective map from $U_\beta$ to the    set of codings:
\begin{equation}\label{eq:U-q}
\mathbf U_\beta:=\set{(d_i)\in \Omega^\N: \sum_{i=1}^\f\frac{d_i}{\beta^i}\in U_\beta}.
\end{equation}
Note that the set-valued map $\beta\mapsto \us_\beta$ is non-decreasing {with respect to the set-inclusion}.  Motivated by the work of periodic unique beta-expansions in $\R$ (see \cite{Allouche-Clarke-Sidorov-2009, Ge-Tan-2017}),  in this paper we will study for each $k\in\N$  the critical  base $\beta_k\in(1,2]$ such that   $\mathbf U_\beta$ contains a   sequence of smallest period $k$ if $\beta>\beta_k$, and $\us_\beta$ contains no   sequence of smallest period $k$ if $\beta<\beta_k$. More precisely, we will determine for each $k\in\N$ the critical base
\begin{equation}\label{eq:beta-k}
\beta_k:=\inf \set{\beta\in(1,2]: \mathbf U_\beta \textrm{ contains a   sequence of smallest period }k}.
\end{equation}
Clearly, for $k=1$ we have $\beta_1=1$, since $\alpha_0^\f=\al_0\al_0\ldots\in\mathbf U_\beta$ for any $\beta\in(1,2]$.

To state our main result  on $\beta_k$ we   define a sequence $(\t_n)$ of blocks in $\set{0,1}^*$ (cf.~\cite[Section 4]{Kong-Li-2020}). Let   $\t_1=100$, and for $n\ge 1$ let
\begin{equation}\label{eq:def-tn}
\t_{n+1}=\t_n^+\Theta(\t_n^+),
\end{equation} where for a block $\c=c_1\ldots c_k\in\set{0,1}^*$ with $c_k=0$ we write $\c^+=c_1\ldots c_{k-1}(c_k+1)$, and the block map $\Theta$ is defined on $\Omega:=\set{000, 001, 100, 101}$ by
\[
 \Theta: \Om\to\Om;\quad 000\mapsto 101, ~001\mapsto 100,~ 100\mapsto 001,~ 101\mapsto 000.
\]
 Then each $\t_n$ is a block of length $3\cdot 2^{n-1}$ by concatenating words from $\Omega$. For example,
 \[
\t_2=101000, \quad\t_3=101001\,000100, \quad \t_4=101001000101\; 000100101000,\]
 and so on. Note by (\ref{eq:def-tn}) that each $\t_{n+1}$ begins with $\t_n^+$. Then the sequence $(\t_n)$ induces a unique componentwise limit
 \begin{equation}\label{eq:lambda-i}
 (\la_i)=\lim_{n\to\f}\t_n=101001000101\ldots\in\set{0,1}^\N.
 \end{equation}
 Furthermore, this limit $(\la_i)$  determines a unique $\beta_c\in(1,2]$ satisfying
\begin{equation}\label{eq:beta-c}
\sum_{i=1}^\f\frac{\la_i}{(\beta_c)^i}=1.
\end{equation}
{We emphasize that $\beta_c\approx 1.55263$ is a transcendental number (see \cite{Kong-Li-2020}).}

  For $\beta\in(1,2]$ we reserve the notation
\[\de(\beta)=\de_1(\beta)\de_2(\beta)\ldots\in\set{0,1}^\N\]
 for the \emph{quasi-greedy} $\beta$-expansion of $1$, which is the lexicographically largest sequence in $\set{0,1}^\N$ not ending with $0^\f$ such that $\sum_{i=1}^\f\frac{\de_i(\beta)}{\beta^i}=1$ (cf.~\cite{Daroczy_Katai_1993}). Some properties of the map $\beta\mapsto\de(\beta)$ are listed in Lemma \ref{lem:delta-beta} below.
It was shown in \cite{Kong-Li-2020} that $(\la_i)$ is indeed the quasi-greedy $\beta_c$-expansion of $1$, i.e., $\de(\beta_c)=(\la_i)$.

\begin{figure}[h!]
  \centering
  \includegraphics[width=10cm]{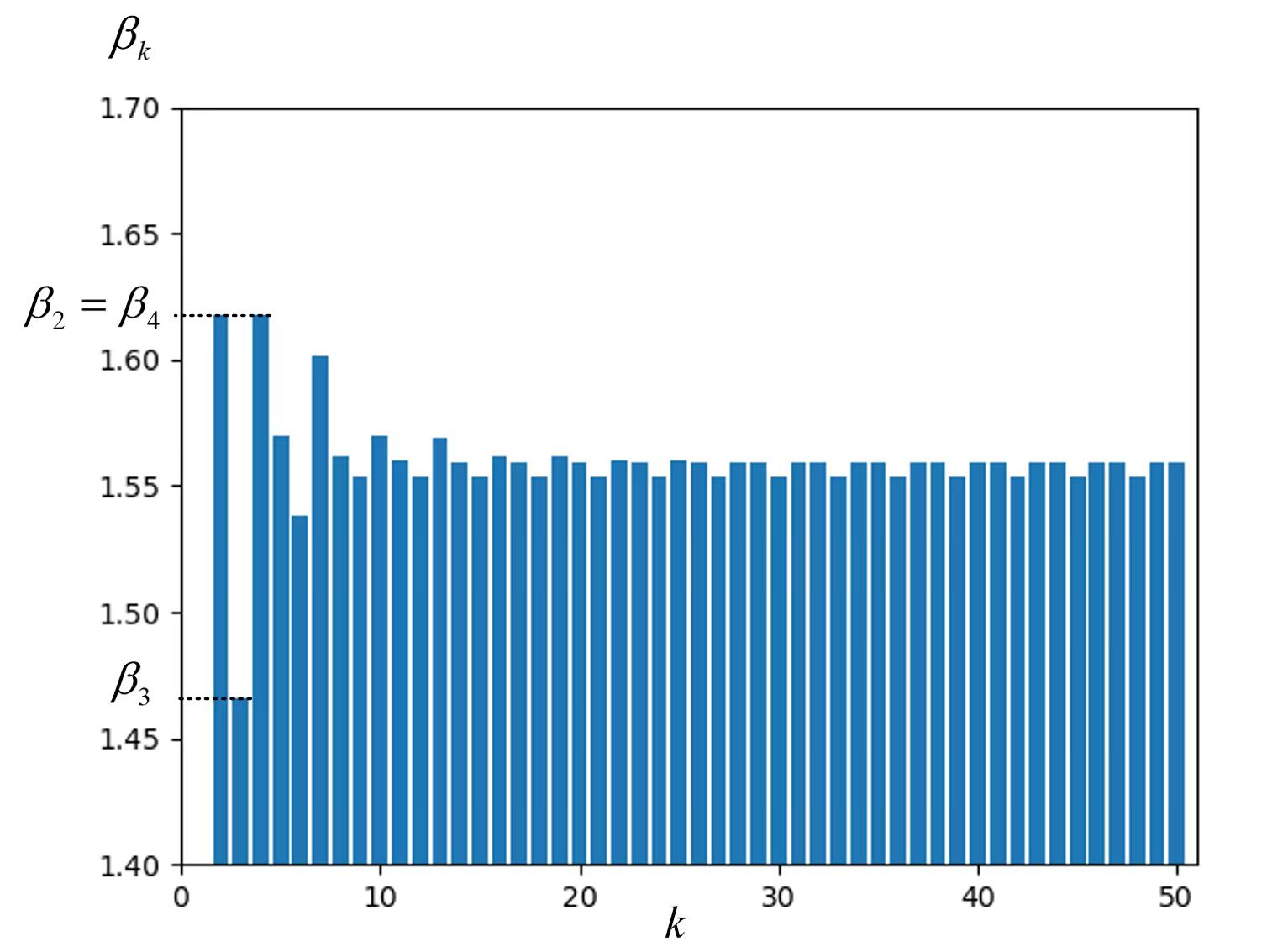}\\
  \caption{The graph of $\beta_k$ for $2\le k\le 50$. Indeed,   $\beta_3\le \beta_k\le\beta_2(=\beta_4)$ for all $k\ge 2$, where $\beta_3\approx 1.46557$ and $\beta_2=\frac{1+\sqrt{5}}{2}\approx 1.61803$.}\label{fig:beta-k}
\end{figure}

Now we state our main result on $\beta_{k}$ for $k\in  \N$  (see Figure \ref{fig:beta-k}). To avoid confusion we write $\N:=\set{1,2,3,\ldots}$ and $\N_0:=\N\cup\set{0}$. {For a real number $r$ let $\lfloor r\rfloor$ be the greatest  integer no larger than $r$, and let $\lceil r\rceil$ be the least integer no smaller than $r$.}
\begin{theorem}
  \label{th:main result}
  The critical base $\beta_k$ for $k\in\N$ is determined as follows.
  \begin{enumerate}
    [{\rm(i)}]

    \item $\beta_1=1$ and $\beta_2=\frac{1+\sqrt{5}}{2}$.

    \item If $k\in 3\N$, then $k=3(2m+1)2^{n}$ for some $m, n\in\N_{0},$ and
  \[
  \de(\beta_{3(2m+1)2^n})=\left\{
  \begin{array}
    {lll}
     \t_{n+1}^\f&\textrm{if}& m=0,\\
     \left(\t_{n+2}^+\Theta(\t_{n+1}^+)\t_{n+2}^{m-1}\right)^\f&\textrm{if}& m\ge 1.
  \end{array}\right.
  \]
  Furthermore, $\beta_{3k}>\beta_{3\ell}$ if and only if $k\rhd \ell$ in the Sharkovskii ordering.

  \item If $k \in 3\N+1$, then $k=3\ell+1$ for some $\ell\in\N$, and
  \[
  \de(\beta_{3\ell+1})=(101(001)^{\lfloor\frac{\ell-1}{2}\rfloor}(010)^{\lceil\frac{\ell-1}{2}\rceil} 0)^\f.
  \]

  \item If $k\in 3\N+2$, then  $k=3\ell+2$ for some $\ell\in\N$, and
  \[
  \de(\beta_{3\ell+2})=(101(001)^{\ell-1}00)^\f.
  \]
  \end{enumerate}

\end{theorem}

\begin{remark}\label{rem:1}\mbox{}
  \begin{enumerate}[{\rm(i)}]
  \item By \cite[Propositions 4.2 and 4.4]{Blanchard-1989} {and Theorem \ref{th:main result}} it follows that each $\beta_k$ is a Perron number.

    \item Note by Theorem \ref{th:main result} (iii) and (iv) that for $\ell=2m+1$ with $m\in\N$ we have
  \[
  \de(\beta_{3(2m+1)+1})=(101(001)^{m}(010)^{m} 0)^\f=(101(001)^{m-1}00)^\f=\de(\beta_{3m+2}).
  \]
  So, \begin{equation}\label{eq:identity}\beta_{6m+4}=\beta_{3m+2}\quad \textrm{for any }m\in\N.\end{equation}
  Indeed, (\ref{eq:identity}) also holds for $m=0$, since $\de(\beta_2)=(10)^\f$ by Theorem \ref{th:main result} (i).

  \end{enumerate}
\end{remark}
{In contrast with the study of periodic unique $\beta$-expansions in $\R$ (see \cite{Allouche-Clarke-Sidorov-2009, Ge-Tan-2017}), the study of periodic unique codings in $S_\beta\subset\R^2$ is more involved. The main difficulty in our proof of Theorem \ref{th:main result} is to show that $\mathbf U_{\beta_k}$ contains no sequence of smallest period $k$. Our strategy is to prove by contradiction based on introducing the notion of admissible blocks (see Definition \ref{def:admissible}). Suppose $\de(\beta_k)=(a_1\ldots a_k)^\f$. Then the block $a_1\ldots a_k$ is   admissible, and it may have multiple representations, which makes our analysis more complicated. Fortunately, when the block $a_1\ldots a_k$ has a prefix $a_1\ldots a_9=\la_1\ldots \la_9=101001000$, we can show that it (nearly) has a unique representation (see Proposition \ref{prop:key-prop}).

The following asymptotic  behavior of the sequence $\set{\beta_k}$ can be deduced from Theorem \ref{th:main result}.}
\begin{proposition}
  \label{prop:asymptotic-beta-k}
  \begin{enumerate}
    [{\rm(i)}]

    \item For any $n\in\N_{0}$ we have
   \[\beta_{3(2m+1)2^n}\searrow \hat\beta_n\quad\textrm{as}\quad \N\ni m\to\f,\]
    where $\de(\hat\beta_n)=\t_{n+2}^+\Theta(\t_{n+2})^\f$. Furthermore, $\hat\beta_n\searrow \beta_c$ as $n\to\f$.

    \item For $m\in\N$ we have
    \[
    \beta_{3(2m+1)2^n}\searrow \beta_c\quad\textrm{as}\quad n\to\f.
    \]Furthermore, for $m=0$ we have $\beta_{3\cdot 2^n}\nearrow\beta_c$ as $n\to\f$.

    \item We have
    \[
    \beta_{3\ell+1}\searrow \beta_a,\quad\beta_{3\ell+2}\searrow\beta_a\quad\textrm{as}\quad\ell\to\f,
    \]
    where $\beta_a\approx 1.55898$ admits $\de(\beta_a)=101(001)^\f$.
  \end{enumerate}
\end{proposition}
Again by \cite[Propositions 4.2 and 4.4]{Blanchard-1989} {and Proposition \ref{prop:asymptotic-beta-k}} one can verify that $\beta_a$ and  $\hat\beta_n$ are Perron numbers {for all $n\in\N$}.
\begin{figure}[h!]
  \centering
  \includegraphics[width=7cm]{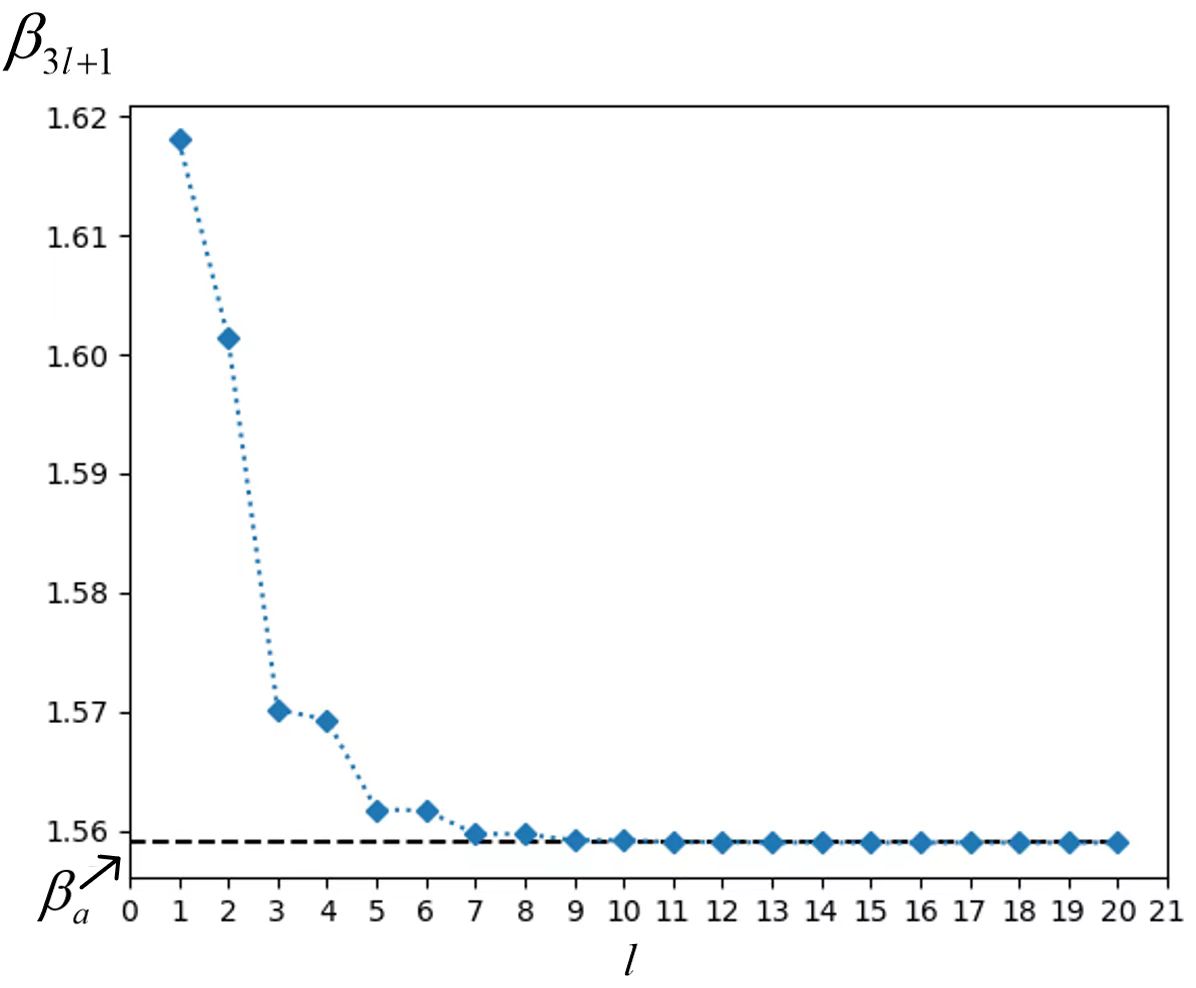}\quad \includegraphics[width=7cm]{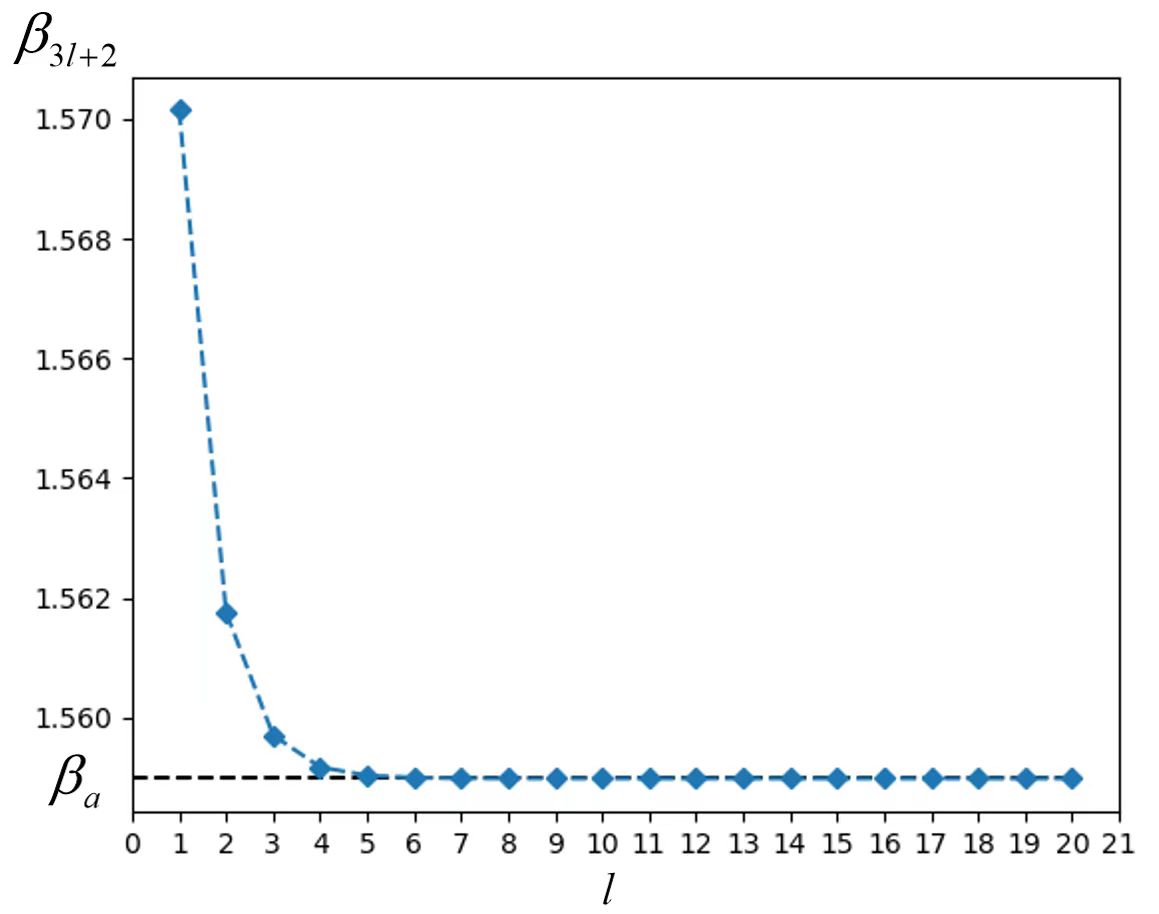}\\
  \caption{Left: the graph of $\beta_{3\ell+1}$ with $1\le \ell\le 20$; right: the graph of $\beta_{3\ell+2}$ with $1\le\ell\le 20$. Indeed, $\beta_{3\ell+1}\searrow\beta_a$, $\beta_{3\ell+2}\searrow\beta_a$ as $\ell\to\f$, where $\beta_a\approx 1.55898$.}\label{fig:beta-3k-1-2}
\end{figure}
 \begin{figure}[h!]
  \centering
  \includegraphics[width=8cm]{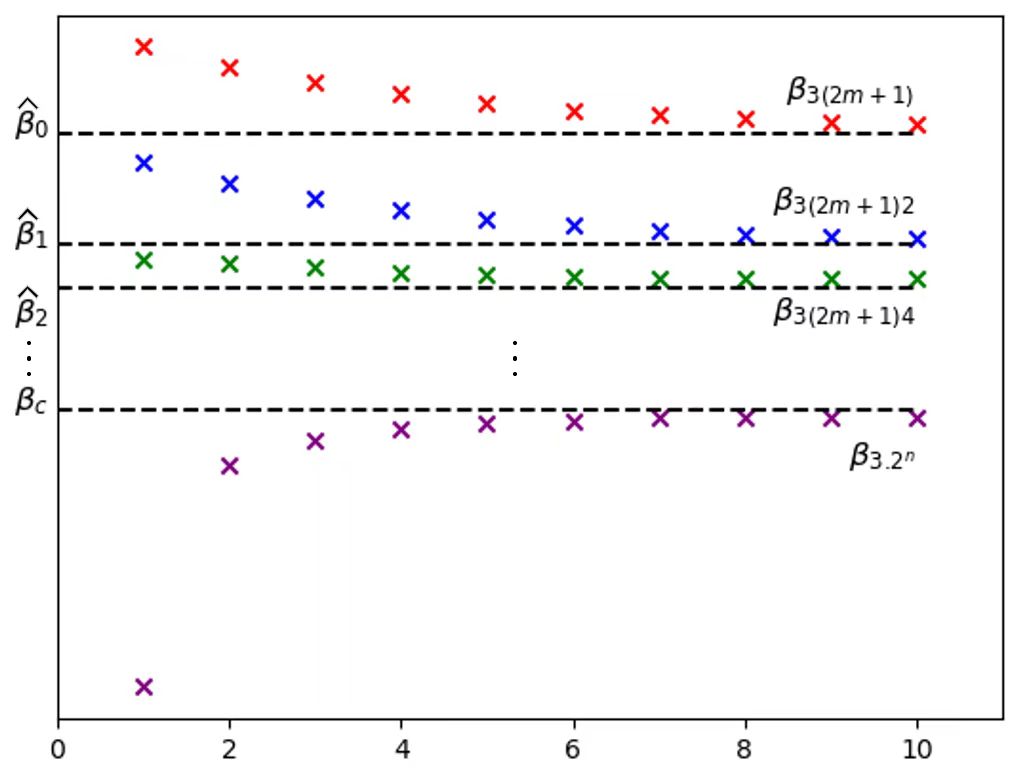}\\
  \caption{The asymptotic behavior  of $\beta_{3(2m+1)2^n}$ with $m, n\in\N_0$. Indeed, $\beta_{3\cdot 2^n}\nearrow\beta_c$ as $n\to\f$; for each $n\in\N_0$, $\beta_{3(2m+1)2^n}\searrow\hat\beta_n$ as $\N\ni m\to\f$; and $\hat\beta_n\searrow\beta_c$ as $n\to\f$.}\label{fig:beta-3k}
\end{figure}

By Proposition \ref{prop:asymptotic-beta-k} (iii) it follows that $\set{\beta_{3\ell+1}}_{\ell=1}^\f$ and $\set{\beta_{3\ell+2}}_{\ell=1}^\f$ are strictly decreasing to the same limit point $\beta_a\approx 1.55898$ (see Figure \ref{fig:beta-3k-1-2}). However, the asymptotic behavior of the sequence $\set{\beta_{3\ell}}_{\ell=1}^\f$ is more interesting (see Figure \ref{fig:beta-3k}). {It has infinitely many accumulation points $\set{\hat\beta_n: n\in\N} \cup\set{\beta_c}$, and $\hat\beta_n\searrow \beta_c$ as $n\to\f$.}
In view of the Sharkoviskii ordering in (\ref{eq:sharkovski}), by Theorem \ref{th:main result} (ii) and Proposition \ref{prop:asymptotic-beta-k} we have the following order {for} $\set{\beta_{3\ell}}_{\ell=1}^\f$:
\begin{equation*}
\begin{array}
  {cccccccccccccc}
      &\beta_{3\cdot3}&>&\beta_{3\cdot 5}&>&\beta_{3\cdot 7}&>&\ldots&>&\beta_{3(2m+1)}&>&\ldots&>&\hat\beta_0\\
  >&\beta_{3\cdot2\cdot 3}&>&\beta_{3\cdot2\cdot 5}&>&\beta_{3\cdot2\cdot 7}&>&\ldots&>&\beta_{3\cdot2(2m+1)}&>&\ldots&>&\hat\beta_1\\
  >&\beta_{3\cdot2^2\cdot 3}&>&\beta_{3\cdot2^2\cdot 5}&>&\beta_{3\cdot2^2\cdot 7}&>&\ldots&>&\beta_{3\cdot2^2(2m+1)}&>&\ldots&>&\hat\beta_{2}\\
  &\vdots&&\vdots&&\vdots&&&&\vdots&&&&\vdots\\
  >&\beta_{3\cdot2^n\cdot 3}&>&\beta_{3\cdot2^n\cdot 5}&>&\beta_{3\cdot2^n\cdot 7}&>&\ldots&>&\beta_{3\cdot2^n(2m+1)}&>&\ldots&>&\hat\beta_n\\
  &\downarrow&&\downarrow&&\downarrow&&&&\downarrow&&&&\downarrow\\
  &\beta_c&&\beta_c&&\beta_c&&&&\beta_c&&&&\beta_c,
\end{array}
\end{equation*}
and
\[
\beta_c>\cdots>\beta_{3\cdot 2^n}>\beta_{3\cdot 2^{n-1}}>\cdots>\beta_{3\cdot 2^2}>\beta_{3\cdot 2}>\beta_3.
\]
\medskip

Observe that $\de(\beta_{3\cdot 3})=(101001000)^\f\prec 101(001)^\f=\de(\beta_a)$. So, by Proposition \ref{prop:asymptotic-beta-k} (iii) it follows that
\begin{align*}
\beta_{3\ell}&<\beta_{3\ell'+1}\le \beta_4\quad\forall \ell, \ell'\in\N;\qquad
 \beta_{3\ell}<\beta_{3\ell'+2}\le \beta_2\quad\forall\ell\in\N, \ell'\in\N_{0}.
 \end{align*}
 Note by Remark \ref{rem:1} that $\beta_2=\beta_4=\frac{1+\sqrt{5}}{2}$. Then $\beta_3\le \beta_k\le \frac{1+\sqrt{5}}{2}$ for any $k\ge 2$. {This means that if $\mathbf U_\beta$ contains a sequence of smallest period $2$ or $4$, then $\mathbf U_\beta$ contains a periodic sequence of smallest period $k$ for any $k\in\N$.}

{It is worth mentioning that the Sierpinski gasket $S_\beta$ studied in our paper is not special. It can be extended to investigate a general class of fat Sierpinski gasket in $\R^2$ generated by homogeneous IFS $\set{g_i(x)=\la x+(1-\la)\mathbf p_i}_{i=1}^3$, where $\la\in(1/2, 1)$ and the vertices $\mathbf p_1, \mathbf p_2, \mathbf p_3$ are not collinear (cf.~\cite{Bro-Mon-Sid-04}).}
The rest of the paper is organized as follows. First we consider in Section \ref{s2} the base $\beta_{3\cdot 2^n}$ for $n\in\N_0$, and show that $\beta_{3\cdot 2^n}\nearrow \beta_c$ as $n\to\f$ (see Theorem \ref{th:lower-higher-orders}). Next we describe in Section \ref{sec:3N} the base $\beta_{3(2m+1)2^n}$ for all $m\in\N$ and $n\in\N_0$. Thus, we have a complete description of $\beta_{3\ell}$ for all $\ell\in\N$, and then establish  Theorem \ref{th:main result} (i)--(ii). In Section \ref{sec:3N+2} we consider $\beta_{3\ell+2}$ for $\ell\in\N$, and prove Theorem \ref{th:main result}   (iv). Finally, in Section \ref{sec:3N+1} we consider $\beta_{3\ell+1}$ for $\ell\in\N$, and prove Theorem \ref{th:main result} (iii). Furthermore, we   prove the asymptotic behavior of $\set{\beta_k}$ as described in Proposition \ref{prop:asymptotic-beta-k}.

 \section{ Periodic sequences in $\mathbf U_\beta$ with $\beta<\beta_c$}\label{s2}
 In this section we will describe the periodic {sequences in $\mathbf U_\beta$} for $\beta<\beta_c$, where $\beta_c$ is defined in (\ref{eq:beta-c}).
First we   need some {terminology} from symbolic dynamics (cf.~\cite{Lind_Marcus_1995}). For a {block} $\c=c_1\ldots c_n\in\set{0,1}^*$ we mean a finite string of zeros and ones. For an integer $k\ge 1$ we denote by $\c^k:=\underbrace{\c\c\cdots\c}_k$ the $k$-times concatenation of $\c$ with itself, and we write for $\c^\f$ the periodic sequence with periodic block $\c$. For a word $\c=c_1\ldots c_n$ with $c_n=0$ we denote by $\c^+:=c_1\ldots c_{n-1}(c_n+1)$.
 For a sequence $(c_i)\in \set{0,1}^\N$ we define its \emph{reflection} by
 $
 \overline{(c_i)}:=(1-c_1)(1-c_2)\cdots.
 $
 {Similarly, for a block $\c=c_1\ldots c_n$ we define its reflection by $\overline{\c}=(1-c_1)(1-c_2)\ldots(1-c_n)$.}
 Throughout the paper we will use the \emph{lexicographical ordering} between sequences and words. More precisely,  for two sequences $(c_i), (d_i)\in\set{0,1}^\N$ we write $(c_i)\prec (d_i)$ or $(d_i)\succ (c_i)$ if $c_1<d_1$, or there exists $k\ge 2$ such that $c_i=d_i$ for all $1\le i<k$ and $c_k<d_k$. Similarly, we write $(c_i)\lle (d_i)$ or $(d_i)\lge(c_i)$ if $(c_i)\prec (d_i)$ or $(c_i)=(d_i)$.  Furthermore, for two words $\c, \d\in\set{0, 1}^*$ we say $\c\prec \d$ if $\c 0^\f\prec \d 0^\f$.

Given $\beta\in(1,2]$, recall that $\de(\beta)=(\de_i(\beta))\in\set{0,1}^\N$ is the {quasi-greedy} $\beta$-expansion of $1$ (cf.~\cite{Daroczy_Katai_1993, Komornik_2011}), which is the lexicographically largest $\beta$-expansion of $1$ not ending with a string of zeros.
 The following characterization of  $\de(\beta)$ was essentially due to Parry \cite{Parry_1960} (see also \cite{ Allaart-Baker-Kong-17, Allaart-2017, DeVries_Komornik_2008}).
\begin{lemma}\label{lem:delta-beta}
\begin{enumerate}
[{\rm(i)}]
\item The map $\beta\mapsto \de(\beta)$ is a strictly increasing bijection from $(1,2]$ onto the set $\mathcal{A}$ of all sequences $(a_i)\in\set{0,1}^\N$ not ending with $0^\f$ and satisfying
\[
a_{n+1}a_{n+2}\ldots \lle a_1 a_2\ldots\quad\forall n\ge 0.
\]

\item The inverse map
\[
\de^{-1}: \mathcal A\ra (1,2];\quad (a_i)\mapsto \de^{-1}((a_i))
\]
is bijective and strictly increasing. Furthermore, $\de^{-1}$ is continuous with respect to the order topology in $\set{0,1}^\N$.
\end{enumerate}
\end{lemma}
 Recall from (\ref{eq:def-tn}) that $\t_n$ is a block of length $3\cdot 2^{n-1}$. Then the sequence  $(a_i)=\t_n^\f$ determines a unique base $\rho_n\in(1,2]$ by
 \begin{equation}\label{eq:rho-n}
\sum_{i=1}^\f\frac{a_i}{(\rho_n)^i}=1.
\end{equation}
 By numerical calculation one can get $\rho_1\approx 1.46557, \rho_2\approx 1.5385, \rho_3\approx 1.55263$, and so on.
   Indeed,   in \cite{Kong-Li-2020} it has shown that $\t_n^\f$ is the  {quasi-greedy} $\rho_n$-expansion of $1$, i.e., $\de(\rho_n)=\t_n^\f$.   For convenience we set $\de(1)=0^\f$. Note that $0^\f\prec\t_1^\f\prec \t_2^\f\prec \t_3^\f\prec\cdots\prec (\la_i)$ in the lexicographical {ordering}, where $(\la_i)=\de(\beta_c)$ is the componentwise limit of $\t_n$.  Therefore, by {(\ref{eq:beta-c}), (\ref{eq:rho-n}) and} Lemma \ref{lem:delta-beta} it follows that
\[
1<\rho_1<\rho_2<\cdots<\rho_n<\rho_{n+1}<\cdots<\beta_c,\quad\textrm{and}\quad \rho_n\nearrow \beta_c.
\]

 In this section we prove the following result.
\begin{theorem}
  \label{th:lower-higher-orders}\mbox{}
Let $\rho_n$ and $\beta_c$ be defined as in (\ref{eq:rho-n}) and (\ref{eq:beta-c}) respectively.
\begin{enumerate}[{\rm(i)}]
  \item For any $n\in\N_0$ we have $\beta_{3\cdot 2^{n}}=\rho_{n+1}$. Thus $\rho_{3\cdot 2^{n}}\nearrow\beta_c$ as $n\to\f$.

\item For any $k\in\N_{\ge 2}\setminus\set{3\cdot 2^{n}: n\in\N_0}$ we have $\beta_k\ge\beta_c$.

  \item $\beta_k\le \beta_2=\frac{1+\sqrt{5}}{2}$ for all $k\in\N$.
\end{enumerate}
\end{theorem}

\begin{remark}\label{rem:2}
Note that $\beta_1=1$. Then by Theorem \ref{th:lower-higher-orders} (i) and (ii) it follows that $\beta_k<\beta_c$ if and only if $k\in\set{1}\cup\set{3\cdot 2^{n}: n\in\N_0}$. Furthermore, by (\ref{eq:beta-k}) and Theorem \ref{th:lower-higher-orders} (iii) we obtain that  for $\beta>\frac{1+\sqrt{5}}{2}$  the   set $\mathbf U_\beta$ contains a sequences of smallest periodic $k$ for any $k\in\N$.
\end{remark}

{In this paper we also need to handle blocks and sequences with alphabet $\Omega=\set{\al_0, \al_1, \al_2}$. Let $\Omega^\N$ be the set of sequences $(d_i)$ with each digit $d_i\in \Omega$. Similarly, for a block $\d=d_1\ldots d_n\in \Omega^*$ and $k\in\N$, we denote by $\d^k=\underbrace{\d\d\ldots \d}_{k}$ the $k$-times concatenation of $\d$. Moreover, we write for $\d^\f\in \Omega^\N$ the periodic sequence with period block $\d$. For a digit $d=(d^1, d^2)\in \Omega=\set{\al_0,\al_1,\al_2}$ we write $d^\+:=d^1+d^2$. Then  $d^1, d^2$ and $d^\+\in\set{0, 1}$. Furthermore,
  \begin{equation}\label{eq:sum-one-digits}
   d^1+d^2+\overline{d^\+}=1\quad \forall d\in \Omega,
   \end{equation}
   where $\overline{d^\+}=1-d^\+$. Note that we are not going to define any order in $\Omega^\N$ or $\Omega^*$.}

Recall that $\U_\beta$ is defined in (\ref{eq:U-q}). Then for each sequence $(d_i)\in \U_\beta$ the point $\sum_{i=1}^\f\frac{d_i}{\beta^i}\in S_\beta$ has a unique coding $(d_i)$.  The following  lexicographical characterization of $\U_\beta$ was given in \cite[Proposition 2.2]{Kong-Li-2020}.

  \begin{lemma}
 \label{lem:lexicographical-characterization-U}
Let $\beta\in(1,2)$. Then     $(d_i)\in\U_\beta$ if and only if   $(d_i)\in \Omega^\N$ satisfies
\[\left\{\begin{array}{lll}
d_{n+1}^1d_{n+2}^1\ldots\prec \de(\beta)&\quad\textrm{whenever}\quad& d_n^1=0,\\
d_{n+1}^2 d_{n+2}^2\ldots\prec \de(\beta)&\quad\textrm{whenever}\quad& d_n^2=0,\\
\overline{d_{n+1}^\+d_{n+2}^\+\cdots} \prec \de(\beta) &\quad\textrm{whenever}\quad& \overline{d_n^\+}=0.
\end{array}\right.\]
 \end{lemma}
 {For $k\in\N$ let ${\Omega(k)}$ be the set of all   sequences in $\Omega^\N$ of smallest period $k$. Based on Lemma \ref{lem:lexicographical-characterization-U} we have a characterization of $\U_\beta\cap {\Omega(k)}$.
 When $k=1$, it is known from \cite{Sidorov_2007} that $\U_\beta\cap \Omega(1)=\set{\al_0^\f, \al_1^\f, \al_2^\f}$ for any $\beta\in(1,2]$. In the following we consider $\U_\beta\cap {\Omega(k)}$ for $k\ge 2$.
 \begin{lemma}
   \label{lem:lexicographical-periodic-U}
   Let $\beta\in(1,2]$ and $(d_i)\in {\Omega(k)}$ with $k\in\N_{\ge 2}$. Then $(d_i)\in\U_\beta$ if and only if
   \[
   \si^n((d_i^1))\prec \de(\beta),\quad \si^n((d_i^2))\prec \de(\beta),\quad \si^n(\overline{(d_i^\+}))\prec \de(\beta)\quad\forall n\ge 0.
   \]
 \end{lemma}
 \begin{proof}
   The sufficiency is obvious by Lemma \ref{lem:lexicographical-characterization-U}. For the necessity we take $(d_i)\in \U_\beta\cap {\Omega(k)}$. By Lemma \ref{lem:lexicographical-characterization-U} it suffices to prove that
   \[
   \si^n((d_i^1))\prec\de(\beta)\quad\textrm{whenever}\quad d_n^1=0\qquad\Longrightarrow\qquad \si^n((d_i^1))\prec \de(\beta)\quad \forall n\ge 0.
   \]
   Since $(d_i)\in {\Omega(k)}$ with $k\ge 2$, by (\ref{eq:sum-one-digits}) it follows that   $(d_i^1)=(d^1_1\ldots d_k^1)^\f\in\set{0,1}^\N$ contains both digits $0$ and $1$. Take $n\in\set{k,k+1,\ldots, 2k-1}$. If   $d_n^1=1$, then there exists a largest $m(<n)$ such that $d_m^1=0$; and thus $d_m^1d_{m+1}^1\ldots d_n^1=01^{n-m}$. By our assumption it follows that
   \begin{align*}
     \si^n((d_i^1))&=d_{n+1}^1d_{n+2}^1\ldots  \lle 1^{n-m}d_{n+1}^1d_{n+2}^1\ldots =\si^m((d_i^1))\prec \de(\beta).
   \end{align*}
   Therefore, $\si^n((d_i^1))\prec \de(\beta)$ for all $k\le n<2k$. Since $(d_i^1)=(d_1^1\ldots d_k^1)^\f$ is periodic,  we conclude that $\si^n((d_i^1))\prec \de(\beta)$ for all $n\ge 0$.

    Similarly, we can prove that $\si^n((d_i^2))\prec \de(\beta)$ and $\si^n(\overline{(d_i^\+)})\prec \de(\beta)$ for all $n\ge 0$. This completes the proof.
 \end{proof}
 }
For $k\in\N$ let
\[B_k:=\set{\beta\in(1,2]: \U_\beta \textrm{ contains a   sequence of smallest period }k}.\]
Based on the lexicographical characterization of $\U_\beta\cap {\Omega(k)}$ we show that each $B_k$ is  a half open interval.
\begin{lemma}\label{lem:basic-form-interval}
  For any $k\in\N$ we have $B_k=(\beta_k,2]$. Furthermore, $\de(\beta_k)=(\de_1(\beta_k)\ldots \de_k(\beta_k))^\f$.
\end{lemma}

\begin{proof}
Note that the set-valued map $\beta\mapsto \U_\beta$ is non-decreasing with respect to the set inclusion. By (\ref{eq:beta-k}) and the definition of $B_k$ it follows that $B_k=(\beta_k,2]$ or $[\beta_k,2]$.
Take $k\in\N$. Note that ${\Omega(k)}$ consists of all periodic sequences in $\Omega^\N$ of {smallest} period $k$. Then $\# {\Omega(k)}\leq (\# \Omega)^k<+\f$. Take $(d_i)=(d_1\ldots d_k)^\f\in {\Omega(k)}$, and    we get three new periodic  sequences: $(d_i^1)$, $(d_i^2)$ and $(d_i^\+)\in\set{0,1}^\N$. Let $(\hat d_i)\in\set{0,1}^\N$ be the lexicographically largest sequence in
\begin{equation*}
\bigcup_{n=0}^{k-1}\set{\si^n((d_i^1)), \si^n((d_i^2)), \si^n(\overline{(d_i^\+)})},
\end{equation*}
where $\sigma$ is the left shift map. {Then by Lemma \ref{lem:lexicographical-periodic-U} it follows that $(d_i)\in\mathbf U_\beta$ if and only if $\de(\beta)\succ(\hat d_i)$.}
Since $\#{\Omega(k)}<+\f$,  the set $\set{ (\hat d_i): (d_i)\in {\Omega(k)}}$ has a   lexicographically smallest sequence, say   $(a_i)\in\set{0,1}^\N$. Then $(a_i)=(a_1\ldots a_k)^\f$. By Lemma \ref{lem:lexicographical-periodic-U} it follows that
\begin{equation}\label{eq:oct-11}
\U_\beta\cap {\Omega(k)}\ne\emptyset\quad \textrm{if and only if}\quad \de(\beta)\succ (a_i).
\end{equation}
 Furthermore, since $(a_i)$ satisfies the inequalities in Lemma \ref{lem:delta-beta} (i),  by (\ref{eq:beta-k}) and (\ref{eq:oct-11}) it follows that $\de(\beta_k)=(a_i)=(a_1\ldots a_k)^\f$. So, by (\ref{eq:oct-11}) and Lemma \ref{lem:delta-beta} it follows that $B_k=(\beta_k,2]$.
\end{proof}
\begin{remark}
  The above proof also provides an algorithm to determine $\beta_k$ for any $k\in\N$. However, when $k$ gets larger, the algorithm does not work efficiently.
\end{remark}
%

\begin{proof}[Proof of Theorem \ref{th:lower-higher-orders}]
First we prove (i).
  Let $\su_0=\al_0, \su_1=\al_1\al_0\al_2$, and for $n\ge 1$ we set
  \[
  \su_{n+1}:=\su_n^{\al_1}\Phi(\su_n^{\al_1}),
  \]
  where $\su_n^{\al_1}$ is the block obtained by changing the last digit of $\su_n$ to $\al_1$, and $\Phi$ is a block map defined on $\Omega=\set{\al_0,\al_1,\al_2}$ by
  \[\Phi: \al_0\mapsto \al_0,\quad \al_1\mapsto \al_2, \quad\al_2\mapsto\al_1.\]
   Then $\su_2=\al_1\al_0\al_1\al_2\al_0\al_2, \su_3=\al_1\al_0\al_1\al_2\al_0\al_1\al_2\al_0\al_2\al_1\al_0\al_2$, and so on. Let $\theta$ be the cyclic permutation on $\Omega$ defined by
   \[\theta: \al_0\mapsto \al_1,\quad \al_1\mapsto \al_2,\quad \al_2\mapsto \al_0.\]
    {Then the map $\theta$ can be applied to blocks by $\theta(d_1\ldots d_m)=\theta(d_1)\ldots\theta(d_m)$.}  {Note that the composition maps} $\theta^2: \al_0\mapsto \al_2, \al_1\mapsto\al_0, \al_2\mapsto \al_1$, and $\theta^3$ is an identity map. By \cite[Proposition 5.5]{Kong-Li-2020} it follows that for any $\beta\in(\rho_n, \rho_{n+1}]$ with $n\in\N$ we have $\U_\beta=\U_{\rho_{n+1}}$, {and} any sequence in $\U_{\rho_{n+1}}$ must end in
  \begin{equation}\label{eq:lower-1}
  \bigcup_{k=0}^n\set{(\su_k)^\f, (\theta(\su_k))^\f, (\theta^2(\su_k))^\f}.
  \end{equation}
  Note that each sequence in $\set{(\su_k)^\f, (\theta(\su_k))^\f, (\theta^2(\su_k))^\f}$ with $1\le k\le n$ is a   sequence of smallest period $3\cdot 2^{k-1}$. Since the set-valued map $\beta\mapsto \U_\beta$ is increasing, it follows that for $\beta\le \rho_{n+1}$ the set $\U_\beta$ can not contain a   sequence of smallest period $3\cdot 2^n$, and thus $\beta_{3\cdot 2^{n}}\ge\rho_{n+1}$.

  On the other hand, by \cite[Lemma 5.3]{Kong-Li-2020} it follows that any sequence in
    \begin{equation*}
  \bigcup_{k=0}^{n+1}\set{(\su_k)^\f, (\theta(\su_k))^\f, (\theta^2(\su_k))^\f}
  \end{equation*}
  belongs to $\U_\beta$ for any $\beta>\rho_{n+1}$. Since $(\su_{n+1})^\f$ is a sequence of smallest period $3\cdot 2^{n}$, we conclude that $ \beta_{3\cdot 2^n}\le \rho_{n+1}$. This proves (i).

  For (ii)  note that $\rho_n\nearrow \beta_c$ as $n\to \f$, and any sequence in $\U_{\rho_{n+1}}$ ends in a periodic sequence in (\ref{eq:lower-1}). Thus, any   sequence of smallest period $k\in\N_{\ge 2}\setminus\set{3\cdot 2^{n}:n\in\N_0}$  can not belong to $\U_\beta$ for $\beta<\beta_c$. In other words, for any $k\in\N_{\ge 2}\setminus\set{3\cdot 2^{n}:n\in\N_0}$   we must have $\beta_k\ge \beta_c$, proving (ii).

Finally we {consider} (iii).
  First we prove $\beta_2\ge \frac{1+\sqrt{5}}{2}$. Suppose on the contrary that $\beta_2<\frac{1+\sqrt{5}}{2}$. Then $\U_{\frac{1+\sqrt{5}}{2}}$ contains a   sequence of smallest period $2$. By symmetry we may assume $(d_i)=(\al_0\al_1)^\f\in\U_{\frac{1+\sqrt{5}}{2}}$. Then
  \[
 d_1^1=0\quad\textrm{and}\quad  d_2^1d_3^1\ldots=(10)^\f=\de(\frac{1+\sqrt{5}}{2}),
  \]
  leading to a contradiction with Lemma \ref{lem:lexicographical-characterization-U}. So,   $\beta_2\ge\frac{1+\sqrt{5}}{2}$ as desired. In the following it suffices to prove $\beta_k\le \frac{1+\sqrt{5}}{2}$ for all $k\in\N$.

  Take $\beta>\frac{1+\sqrt{5}}{2}$. Then $\de(\beta)\succ (10)^\f$. We only need to prove that for any $k\in\N$ the set $\U_\beta$ contains a  sequence of smallest period $k$. By Lemma \ref{lem:lexicographical-characterization-U} one can   verify that $\al_0^\f, (\al_0\al_1)^\f\in\U_\beta$, and then $\beta_1, \beta_2\le \beta$. For $k\ge 3$ we can write $k=3\ell+r$ with $\ell\in\N$ and $r\in\set{0,1,2}$. {We} split the proof of $\beta_k\le \beta$ into the following three cases.

  (I). $k=3\ell$ with $\ell\in\N$. Let \[(d_i)=((\al_0\al_1\al_2)^{\ell-1}\al_1\al_0\al_2)^\f.\] Then $(d_i)$ is a   sequence of smallest period $3\ell$. Furthermore, $(d_i^1)=((010)^{\ell-1}100)^\f$ and thus $\si^n((d_i^1))\prec (10)^\f\prec \de(\beta)$ for any $n\ge 0$. Similarly, $\si^n((d_i^2))=\si^n((001)^\f)\prec \de(\beta)$ for all $n\ge 0$, and $\si^n(\overline{(d_i^{\+})})=\si^n((100)^{\ell-1}010)^\f\prec {\de(\beta)}$ for all $n\ge 0$. By Lemma \ref{lem:lexicographical-characterization-U} it follows that $(d_i)\in \U_\beta$, and thus $\beta_{3\ell}\le\beta$.

  (II). $k=3\ell+1$ with $\ell\in\N$. Let \[(d_i)=((\al_0\al_1\al_2)^\ell \al_1)^\f.\] Then $(d_i)$ is a  sequence of smallest period $3\ell+1$. Note that $(d_i^1)=((010)^\ell 1)^\f, (d_i^2)=((001)0)^\f$ and $\overline{(d_i^{\+})}=((100)^\ell 0)^\f$. By Lemma \ref{lem:lexicographical-characterization-U} and {using} $\de(\beta)\succ (10)^\f$ one can verify that $(d_i)\in\U_\beta$, and thus $\beta_{3\ell+1}\le\beta$.

  (III). $k=3\ell+2$. Let \[(d_i)=((\al_0\al_1\al_2)^\ell \al_0\al_1)^\f.\] Then $(d_i)$ is a  sequence of smallest period $3\ell+2$. Note that $(d_i^1)=((010)^\ell 01)^\f, (d_i^2)=((001)^\ell 00)^\f$ and $\overline{(d_i^{\+})}=((100)^\ell 10)^\f$. {Again,} by Lemma \ref{lem:lexicographical-characterization-U} and {using} $\de(\beta)\succ (10)^\f$ it follows that $(d_i)\in\U_\beta$, and hence $\beta_{3\ell+2}\le\beta$.

  By Cases (I)--(III)  and using that $\beta_1, \beta_2\le \beta$ it follows that $\beta_k\le \beta$ {for all $k\in\N$}. Since $\beta>\frac{1+\sqrt{5}}{2}$ was arbitrary, we conclude that $\beta_k\le \frac{1+\sqrt{5}}{2}$, completing the proof.
\end{proof}

\section{Minimal base for $\U_\beta$ to have a sequence of period  $3\N$}\label{sec:3N}
In this section we will determine $\beta_k$ for all $k\in 3\N$. Note by (\ref{eq:def-tn}) and (\ref{eq:lambda-i}) that
$
\la_1\ldots \la_{3\cdot 2^n}^-=\t_{n+1}$ for all $n\in\N_0.
$
Then by {(\ref{eq:rho-n}) and} Theorem \ref{th:lower-higher-orders} it follows that
\[
\de(\beta_{3\cdot 2^n})=(\la_1\ldots \la_{3\cdot 2^n}^-)^\f{=\t_{n+1}^\f}.
\]
 For other $k\in3\N$ we have the following result.
\begin{theorem}
  \label{th:eta-period-3n}
  For any $m\in\N$ and $n\in\N_0$ we have
  \[
  \de(\beta_{3(2m+1)2^n})=\left(\la_1\ldots\la_{3(2^{n+1}+2^n)}(\la_1\ldots \la_{3\cdot 2^{n+1}}^-)^{m-1}\right)^\f.
  \]
\end{theorem}
\begin{remark}\mbox{}
\begin{enumerate}[{\rm(i)}]
\item
  It is worth mentioning  that $\beta_{3(2m+1)2^n}>\beta_c$ for any $m\in\N, n\in\N_0$. This can be deduced by Lemma \ref{lem:delta-beta} and the following observation:
  \begin{align*}
    \de(\beta_{3(2m+1)2^n})&=\left(\la_1\ldots\la_{3(2^{n+1}+2^n)}(\la_1\ldots \la_{3\cdot 2^{n+1}}^-)^{m-1}\right)^\f\\
    &=\left(\la_1\ldots\la_{3\cdot 2^n}\Theta(\la_1\ldots \la_{3\cdot 2^n})^+\Theta(\la_1\ldots \la_{3\cdot 2^n})\Big(\la_1\ldots \la_{3\cdot 2^n}\Theta(\la_1\ldots \la_{3\cdot 2^n})\Big)^{m-1}\right)^\f\\
    &\succ \la_1\ldots\la_{3\cdot 2^n}\Theta(\la_1\ldots \la_{3\cdot 2^n})^+\Theta(\la_1\ldots \la_{3\cdot 2^n})\la_1\ldots\la_{3\cdot 2^n}\Theta(\la_1\ldots \la_{3\cdot 2^n})\la_1\ldots \la_{3\cdot 2^n}^- 1^\f\\
    &\succ(\la_i)=\de(\beta_c).
  \end{align*}

  \item Note that $\t_{n+1}^+=\la_1\ldots\la_{3\cdot 2^n}$ for all $n\in\N_0$. Then by (\ref{eq:def-tn})  we can rewrite Theorem \ref{th:eta-period-3n}   as
  \[
  \de(\beta_{3(2m+1)2^n})=\left(\t_{n+2}^+\Theta(\t_{n+1}^+)\t_{n+2}^{m-1}\right)^\f,
  \]
  which coincides with Theorem \ref{th:main result} (ii).
  \end{enumerate}
\end{remark}

For $n\in\N_0$ we write $\ga_1\ldots \ga_{3\cdot 2^n}:=\Theta(\la_1\ldots \la_{3\cdot 2^n})$. Then
\begin{equation}\label{eq:la-ga}
\la_{3i+1}+\ga_{3i+1}=1,\quad \la_{3i+2}+\ga_{3i+2}=0\quad\textrm{and}\quad \la_{3i+3}+\ga_{3i+3}=1\quad\forall 0\le i<2^n.
\end{equation} The following lemma was proven in \cite[Lemma 4.1]{Kong-Li-2020}.
\begin{lemma}
  \label{lem:inequalities}
  For any $n\in\N_0$ we have
  \begin{align*}
    \ga_1\ldots \ga_{3\cdot 2^n-i}&\prec \la_{i+1}\ldots \la_{3\cdot 2^n}\lle \la_1\ldots \la_{3\cdot 2^n-i},\\
    \ga_1\ldots \ga_{3\cdot 2^n-i}&\lle \ga_{i+1}\ldots \ga_{3\cdot 2^n}\prec \la_1\ldots \la_{3\cdot 2^n-i}
  \end{align*}
  for all $0\le i<3\cdot 2^n$.
\end{lemma}
First we prove Theorem \ref{th:eta-period-3n} for $m=1$.

\begin{proposition}\label{prop:m=1}
For any    $n\in\N_0$ we have $\de(\beta_{9\cdot 2^n})=(\la_1\ldots \la_{9\cdot 2^n})^\f$.
\end{proposition}

\begin{proof}Let $k=9\cdot 2^n$ with $n\in\N_0$. Note by Lemma \ref{lem:basic-form-interval} that $\de(\beta_k)$ is a periodic sequence of period $k$. Then
  by Theorem \ref{th:lower-higher-orders} (ii) and Lemma \ref{lem:delta-beta} it follows that
  \[\de(\beta_k)\lge (\la_1\ldots \la_k)^\f.\]
   So in the following it suffices to prove that for any $\beta\in(1,2]$ satisfying $\de(\beta)\succ (\la_1\ldots \la_k)^\f$, the set $\U_\beta$ contains  a sequence of smallest period $k$.

  Take $\beta\in(1,2]$ so that $\de(\beta)\succ (\la_1\ldots \la_k)^\f$. Let $(d_i)=(d_1\ldots d_k)^\f$ such that
  \begin{equation}\label{eq:feb15-1}
  (d_i^1)=(\la_1\ldots \la_k)^\f,\quad (d_i^2)=(\ga_1\ldots \ga_k)^\f,\quad (\overline{d_i^\+})=(010)^\f.
  \end{equation}Then by (\ref{eq:la-ga}) it follows that $d_i^1+d_i^2+\overline{d_i^\+}=1$ for all $i\in\N$. Thus, $(d_i)\in \set{\al_0, \al_1,\al_2}^\N$.
  Furthermore, by the definition of $\la_1\ldots \la_k$ one can verify that $(d_i)$ has smallest period $k$. So, $(d_i)\in {\Omega(k)}$. In the following it suffices to prove  that $(d_i)\in\us_\beta$.
  Clearly,  since $\la_1\la_2\la_3=101$, we have  $\si^n((\overline{d_i^\+}))=\si^n((010)^\f)\prec (\la_1\ldots \la_k)^\f\prec\de(\beta)$. Therefore,  by Lemma \ref{lem:lexicographical-periodic-U} and (\ref{eq:feb15-1})  we only need  to prove that
\begin{equation}
    \label{eq:di-Ubeta-1}
    \si^j((\la_1\ldots \la_k)^\f)\prec \de(\beta)\quad\forall 0\le j<k,
    \end{equation}
    and
    \begin{equation}
    \label{eq:di-Ubeta-2}
    \si^j((\ga_1\ldots \ga_k)^\f)\prec \de(\beta)\quad\forall 0\le j<k.
    \end{equation}

 First we prove (\ref{eq:di-Ubeta-1}).
    Clearly, (\ref{eq:di-Ubeta-1}) holds for $j=0$ by our assumption.
     Note that
  \begin{equation}\label{eq:di-1-2}
  \la_1\ldots \la_k=\la_1\ldots \la_{3\cdot 2^n}\ga_1\ldots \ga_{3\cdot 2^n}^+\ga_1\ldots\ga_{3\cdot 2^n}.
  \end{equation}
     Since $\ga_1=0<\la_1$,  (\ref{eq:di-Ubeta-1})   holds for $j=3\cdot 2^n$ and $j=3\cdot 2^{n+1}$. Furthermore, by   Lemma \ref{lem:inequalities} it follows that for any $1\le j<3\cdot 2^n$,
     \begin{align*}
     \la_{j+1}\ldots \la_{3\cdot 2^n}\ga_1\ldots \ga_j&\prec \la_1\ldots \la_{3\cdot 2^n},\\
      \ga_{j+1}\ldots \ga_{3\cdot 2^n}^+\ga_1\ldots\ga_j&\prec \la_1\ldots \la_{3\cdot 2^n},\\
      \ga_{j+1}\ldots\ga_{3\cdot 2^n}&\prec\la_1\ldots \la_{3\cdot 2^n-j}.
     \end{align*}
    This, together with (\ref{eq:di-1-2}), proves (\ref{eq:di-Ubeta-1}) for all $0\le j<k$.

     Next we prove (\ref{eq:di-Ubeta-2}).
   Note that $\ga_1\ldots \ga_k=\ga_1\ldots \ga_{3\cdot 2^{n+1}}\la_1\ldots \la_{3\cdot 2^n}$. Then (\ref{eq:di-Ubeta-2}) follows
 by  Lemma \ref{lem:inequalities}   that
 \begin{align*}
  \ga_{j+1}\ldots \ga_{3\cdot 2^{n+1}}&\prec \la_1\ldots \la_{3\cdot 2^{n+1}-j}\quad\forall 0\le j<3\cdot 2^{n+1}, \\
 \quad \la_1\ldots \la_{3\cdot 2^n}\ga_1\ldots \ga_{3\cdot 2^n}&= \la_1\ldots \la_{3\cdot 2^n}^-\prec \la_1\ldots \la_{3\cdot 2^n},
 \end{align*}
 and
 \[
 \la_{j+1}\ldots\la_{3\cdot 2^n}\ga_1\ldots \ga_{j}\prec \la_1\ldots \la_{3\cdot 2^n}\quad \forall 1\le j<3\cdot 2^n.
 \]
 This completes the proof.
 \end{proof}
In the following we determine $\beta_k$ for $k=3(2m+1)2^n$ with $m\in\N_{\ge 2}$ and $n\in\N_0$.
The upper bound for $\beta_k$ can be obtained similar to the proof of Proposition \ref{prop:m=1}.

\begin{lemma}
  \label{lem:upper-bound-eta-3k}
  If $k=3(2m+1)2^n$ with $m\in\N_{\ge 2}$ and $n\in\N_0$, then
  \[
  \de(\beta_k)\lle\left(\la_1\ldots\la_{3(2^{n+1}+2^n)}(\la_1\ldots \la_{3\cdot 2^{n+1}}^-)^{m-1}\right)^\f.
  \]
\end{lemma}
\begin{proof}
 Let $\beta\in(1,2]$ so that
  \begin{equation}\label{eq:feb15-2}
  \de(\beta)\succ \left(\la_1\ldots\la_{3(2^{n+1}+2^n)}(\la_1\ldots \la_{3\cdot 2^{n+1}}^-)^{m-1}\right)^\f.
  \end{equation} It suffices to prove that $\U_\beta$ contains a sequence of smallest period $k=3(2m+1)2^n$.
  Let $(d_i)=(d_1\ldots d_k)^\f$ such that
  \begin{align*}
    (d_i^1)&=\left(\la_1\ldots\la_{3(2^{n+1}+2^n)}(\la_1\ldots \la_{3\cdot 2^{n+1}}^-)^{m-1}\right)^\f,\\
    (d_i^2)&=\left(\ga_1\ldots\ga_{3(2^{n+1}+2^n)}(\ga_1\ldots \ga_{3\cdot 2^{n+1}}^+)^{m-1}\right)^\f \quad\textrm{and}\quad
    (\overline{d_i^\+})=(010)^\f.
  \end{align*}
  Then one can verify that $(d_i)\in {\Omega(k)}$. So in the following we only need to prove $(d_i)\in \U_\beta$.
   Note that $\la_1\la_2\la_3=101$. By (\ref{eq:feb15-2}) it is clear that $\si^j((\overline{d_i^\+}))\prec \de(\beta)$ for all $j\ge 0$. Therefore, by Lemma \ref{lem:lexicographical-characterization-U} it suffices to prove that $\si^j((d_i^1))\prec \de(\beta)$ and $\si^j((d_i^2))\prec \de(\beta)$ for all $j\ge 0$.

  First we prove $\si^j((d_i^1))\prec \de(\beta)$ for all $j\ge 0$. By (\ref{eq:feb15-2}) it is clear that $\si^j((d_i^1))\prec\de(\beta)$ for $j=0$. For $j\ge 1$  note that
   \begin{equation}\label{eq:august-11-1}
   (d_i^1) =\left(\la_1\ldots \la_{3\cdot 2^{n+1}}(\ga_1\ldots \ga_{3\cdot 2^{n+1}}^+)^{m-1}\ga_1\ldots \ga_{3\cdot 2^n}\right)^\f.
   \end{equation}
     By Lemma \ref{lem:inequalities} it follows that
   \[\la_{i+1}\ldots \la_{3\cdot 2^{n+1}}\ga_1\ldots \ga_i\prec \la_1\ldots\la_{3\cdot 2^{n+1}}\quad\forall 1\le i<3\cdot 2^{n+1}.\]
   This, together with (\ref{eq:feb15-2}) and (\ref{eq:august-11-1}), implies that  $\si^j((d_i^1))\prec \de(\beta)$ for all $0\le j<3\cdot 2^{n+1}$. Furthermore, by Lemma \ref{lem:inequalities} we have
   \begin{equation*}
     \label{eq:july15-1}
     \begin{split}
     \ga_{i+1}\ldots \ga_{3\cdot 2^{n+1}}^+\ga_1\ldots \ga_{i}&\prec \la_1\ldots \la_{3\cdot 2^{n+1}}\quad\forall 0\le i<3\cdot 2^{n+1},\\
     \ga_{i+1}\ldots \ga_{3\cdot 2^n}&\prec \la_1\ldots\la_{3\cdot 2^n-i}\quad \forall 0\le i<3\cdot 2^n.
     \end{split}
   \end{equation*}
   Then by (\ref{eq:feb15-2}) and (\ref{eq:august-11-1}) it follows that $\si^j((d_i^1))\prec \de(\beta)$ for all $3\cdot 2^{n+1}\le j\le 3(3^{n+1}m+2^n)$.
    Since $(d_i^1)$ is periodic, this proves $\si^j((d_i^1))\prec \de(\beta)$ for all $j\ge 0$.

  Next we prove $\si^j((d_i^2))\prec \de(\beta)$ for all $j\ge 0$. Observe that
   \begin{align*}
  (d_i^2)=\left(\ga_1\ldots\ga_{3\cdot 2^{n+1}}(\la_1\ldots \la_{3\cdot 2^{n+1}}^-)^{m-1}\la_1\ldots \la_{3\cdot 2^n}\right)^\f
   \end{align*}
   begins with $\ga_1\ldots \ga_{3\cdot 2^{n+1}}(\la_1\ldots \la_{3\cdot 2^{n+1}}^-)^m$.
   Then by (\ref{eq:feb15-2}), Lemma \ref{lem:inequalities} and the same argument as above one can verify that $\si^j((d_i^2))\prec \de(\beta)$ for all $j\ge 0$.
   \end{proof}

The proof for the lower bound of $\beta_k$ with $k=3(2m+1)2^n$ is more involved, where $m\in\N_{\ge 2}$ and $n\in\N_0$. {For this we introduce the notion of admissible blocks.}
\begin{definition}\label{def:admissible}
  A   block $a_1\ldots a_k\in\set{0,1}^*$ is called \emph{admissible} if there exists an aperiodic  block $d_1\ldots d_k\in \Omega^k$ satisfying
  \[
  \overline{d_1^\+ \ldots d_k^\+}\lle d_1^2\ldots d_k^2 \lle d_1^1\ldots d_k^1=a_1\ldots a_k
  \]
  and
  \begin{equation}\label{eq:admissible-condition}
  \begin{split}
d_{j+1}^1\ldots d_k^1d_1^1\ldots d_j^1&\lle a_1\ldots a_{k}\quad \forall 0\le j<k,\\
        d_{j+1}^2\ldots d_k^2d_1^2\ldots d_j^2&\lle a_1\ldots a_{k}\quad  \forall 0\le j<k,\\
     \overline{d_{j+1}^\+\ldots d_k^\+ d_1^\+\ldots d_j^\+} &\lle a_1\ldots a_{k}\quad  \forall 0\le j<k.  \end{split}
  \end{equation}
\end{definition}Note by Lemma \ref{lem:delta-beta} that the first line of inequalities in (\ref{eq:admissible-condition}) guarantee that $(a_1\ldots a_k)^\f$  is the quasi-greedy expansion of $1$ for some base $\beta\in(1,2]$.
In Definition \ref{def:admissible} we call the block $d_1\ldots d_k\in \Omega^k$ a \emph{representation} of the admissible block $a_1\ldots a_k\in\set{0,1}^k$. In general, an admissible block $a_1\ldots a_k$ may have multiple representations (see Example \ref{ex:1}).
However, if the admissible block $a_1\ldots a_k$   has a prefix $a_1\ldots a_{9}=\la_1\ldots \la_{9}=101001000$, then it has a unique representation (see Proposition \ref{prop:key-prop} below).
\begin{example}
  \label{ex:1}
  Let $\a=a_1\ldots a_9=110101000$, and take
  \[
  \c=\al_1\al_1\al_2\al_1\al_0\al_1\al_2\al_0\al_2\quad\textrm{and}\quad \d=\al_1\al_1\al_2\al_1\al_2\al_1\al_0\al_2\al_0.
  \]
  Then
  \begin{align*}
  \c^1&=110101000=\a,\quad \c^2=001000101,\quad\overline{\c^\+}=000010010;\\
    \d^1&=110101000=\a,\quad \d^2=001010010,\quad\overline{\d^\+}=000000101.
     \end{align*}
One can verify that both $\c$ and $\d$ are representations of $\a$.  
\end{example}
Observe by  the proof of Lemma \ref{lem:basic-form-interval} that if $\delta(\beta_k)=(a_1\ldots a_k)^\f$, then $a_1\ldots a_k$ must be an admissible block.
 So to prove the lower bound of $\beta_k$ in Theorem \ref{th:eta-period-3n} it suffices to prove that for $k=3(2m+1)2^n$ with $m\in\N_{\ge 2}$ and $n\in\N_0$ there is no admissible block $a_1\ldots a_k$ satisfying
\begin{equation}\label{eq:etak-lower-bound}
\de(\beta_c) \lle (a_1\ldots a_k)^\f \prec(\la_1\ldots \la_{3(2^n+2^n)}(\la_1\ldots \la_{3\cdot 2^n}^-)^{m-1})^\f.
\end{equation}
{In view of Example \ref{ex:1}, the multiple representations of an admissible block complicates our proof of the lower bound of $\beta_k$.}

\begin{figure}[h!]
  \centering
  \begin{tikzpicture}[->,>=stealth',shorten >=1pt,auto,node distance=3cm, semithick,scale=3]

  \tikzstyle{every state}=[minimum size=0pt,fill=none,draw=black,text=black]

  \node[state] (A)                    {$001$ };
  \node[state]         (B) [ right of=A] { $000$};
  \node[state]         (C) [ above of=A] {$101$};
  \node[state](D)[right of=C]{$100$};

  \path[->,every loop/.style={min distance=0mm, looseness=25}]
  (C) edge[->,left] node{} (A)
  (C) edge[bend right,->,right] node{} (B)
(D) edge[->,above] node{} (C)
(A)edge[->,below] node{} (B)
(B) edge[bend right,->,right] node{} (C)
(B) edge[->,right] node{} (D)
;
\end{tikzpicture}
\caption{The directed graph representing the subshift of finite type $X$.}
\label{fig:directed-graph}
\end{figure}
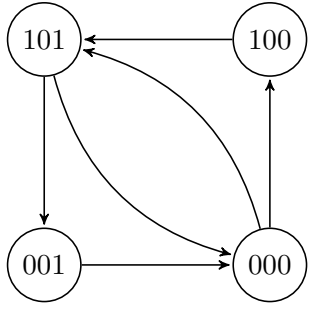

Let $X$ be the subshift of finite type over the states $\set{000, 001, 100, 101}$
represented by the directed graph   in Figure \ref{fig:directed-graph}.
One can verify  that the sequence $\de(\beta_c)=(\la_i)=101001000101\ldots$ is in $X$.
Let $\mathcal B^*(X)$ be the set of all possible blocks appearing in sequences of $X$. The following   key proposition, which plays an essential role in our proof of Theorem \ref{th:eta-period-3n}, shows that any admissible block beginning with $\la_1\ldots \la_{9}$ must belong to $\mathcal B^*(X)$, {and it has a unique representation}.
\begin{proposition}[Key proposition]
  \label{prop:key-prop}
  {If $a_1\ldots a_{k}$ is an admissible block with $k\ge 9$, and it has a prefix
  $
  a_1\ldots a_{9}=\la_1\ldots \la_{9}=101001000,
 $
  then}
  \[a_1\ldots a_{3\ell} \in \mathcal B^*(X),\]
  where    $  \ell= \lfloor k/3\rfloor$. Furthermore, for any representation $d_1\ldots d_k$ of $a_1\ldots a_k$ we have
\[d_1^1\ldots d_{3\ell}^1=a_1\ldots a_{3\ell},\quad d^2_1\ldots d^2_{3\ell}=(010)^\ell\quad\textrm{and}\quad \overline{d_1^\+ \ldots  d_{3\ell}^\+}=\Theta(a_1\ldots a_{3\ell}).\]
\end{proposition}
\begin{proof}
{Since $a_1\ldots a_{k}$ is admissible, it has a representation $d_1\ldots d_{k}\in \Omega^{k}$. Then $d_1^1\ldots d_{k}^1=a_1\ldots a_{k}$ and $d_1\ldots d_{k}$ satisfies (\ref{eq:admissible-condition}). Let $\ell= \lfloor k/3\rfloor$. Then $\ell\ge 3$, and $k\in\set{3\ell, 3\ell+1, 3\ell+2}$.} By (\ref{eq:sum-one-digits}) it suffices to prove $a_1\ldots a_{3\ell}\in \mathcal B^*(X)$ and $d_1^2\ldots d_{3\ell}^2=(010)^\ell$.

First we prove $a_1\ldots a_9\in \mathcal B^*(X)$ and $d_1^2\ldots d_9^2=(010)^3$. Clearly, $a_1\ldots a_9=101001000\in \mathcal B^*(X)$. {Then it suffices to prove $d_1^2\ldots d_9^2=(010)^3$.} Note that
 \begin{equation*}\label{eq:key-1}
  d_1^1\ldots d_{9}^1=\la_1\ldots \la_{9}\lge d_1^2\ldots d_9^2\lge \overline{d_1^\+}\ldots \overline{d_9^\+}.
  \end{equation*}
Since $d_1^1d_2^1d_3^1=101$, by {(\ref{eq:admissible-condition})} and (\ref{eq:sum-one-digits}) it follows that $d_1^2d_2^2d_3^2=010$  and   $\overline{d_1^\+d_2^\+d_3^\+}=000$. Furthermore, by using $a_1\ldots a_9=101001000$ and (\ref{eq:admissible-condition}) {we obtain} that
\begin{equation}\label{eq:forbidden-blocks}
\begin{split}
&\textrm{ the blocks }\quad 11, \quad 10101,\quad 10100101\quad\textrm{ and }\quad 101001001\\
&\textrm{ are all forbidden in } \quad d_1^1\ldots d_{k}^1, \quad d_1^2\ldots d_{k}^2\quad \textrm{ and }\quad\overline{d_1^\+\ldots d_{k}^\+}.
\end{split}
 \end{equation}
Since $a_4a_5a_6=001$, by (\ref{eq:forbidden-blocks}) and (\ref{eq:sum-one-digits}) it follows that $d_4^2d_5^2d_6^2\in\set{010, 100}$. {Note by (\ref{eq:forbidden-blocks}) that $11$ is forbidden in $d_1^2\ldots d_9^2$ and $\overline{d_1^\+\ldots d_9^\+}$.}  If $d_4^2d_5^2d_6^2=100$, then by {(\ref{eq:sum-one-digits})} and $a_7a_8a_9=000$ it follows that
 \[
 \left(\begin{array}{c}
   d_1^1\ldots d_{9}^1\\
   d_1^2\ldots d_9^2\\
   \overline{d_1^\+ \ldots d_9^\+}
 \end{array}\right)=\left(
 \begin{array}{c}
   101\;001\;000\\
   010\;100\;101\\
   000\;010\;010
 \end{array}\right)\quad\textrm{or}\quad
  \left(\begin{array}{c}
   d_1^1\ldots d_{9}^1\\
   d_1^2\ldots d_9^2\\
   \overline{d_1^\+ \ldots d_9^\+}
 \end{array}\right)=\left(
 \begin{array}{c}
   101\;001\;000\\
   010\;100\;010\\
   000\;010\;101
 \end{array}\right).
 \]
 While in the first case we have $d_2^2\ldots d_9^2=10100101$, and in the second case we have $\overline{d_5^\+}\ldots \overline{d_9^\+}=10101$. Both  cases lead to a contradiction with {(\ref{eq:forbidden-blocks})}. Therefore,
 $d_4^2d_5^2 d_6^2=010$. Note that $a_7a_8a_9=000$ {and $11$ is forbidden by (\ref{eq:forbidden-blocks})}.  By (\ref{eq:forbidden-blocks}) and (\ref{eq:sum-one-digits}) we obtain that
 \[
 \left(\begin{array}{c}
   d_1^1\ldots d_{9}^1\\
   d_1^2\ldots d_{9}^2\\
   \overline{d_1^\+ \ldots d_{9}^\+}
 \end{array}\right)=\left(
 \begin{array}{ccc }
   101&001&000 \\
   010&010&101 \\
   000&100&010
 \end{array}\right)\quad\textrm{or}\quad  \left(\begin{array}{c}
   d_1^1\ldots d_{9}^1\\
   d_1^2\ldots d_{9}^2\\
   \overline{d_1^\+ \ldots d_{9}^\+}
 \end{array}\right)=\left(
 \begin{array}{ccc }
   101&001&000 \\
   010&010&010 \\
   000&100&101
 \end{array}\right).
 \]
  If   $d_7^2d_8^2d_9^2 =101$, then we will have $d_5^2\ldots d_9^2=10101$, {again} leading to a contradiction {with (\ref{eq:forbidden-blocks})}. So   $d_7^2 d_8^2 d_9^2=010$ {as desired}.

{If $\ell=3$, then we are done by the above arguments. In the following we assume $\ell>3$. By the induction process we
assume} $a_1\ldots a_{3m}\in \mathcal B^*(X)$ for some $3\le m<\ell$ and $d_1^2\ldots d_{3m}^2=(010)^m$, {and we are going to} prove $a_1\ldots a_{3m+3}\in\mathcal B^*(X)$ and $d_{3m+1}^2d_{3m+2}^2d_{3m+3}^2=010$.

Since $a_1\ldots a_{3m}\in\mathcal B^*(X)$, it follows that either $a_{3m-8}\ldots a_{3m}$ or $\Theta(a_{3m-8}\ldots a_{3m})$ belongs to the set
  \begin{equation}\label{eq:H}
  H:=\set{101001000, 101000101, 101000100, 100101000, 100101001}.
  \end{equation}
 By symmetry we {may assume}   $a_{3m-8}\ldots a_{3m}\in H$, {and the proof will be split into the following five cases}.

 Case I. $a_{3m-8}\ldots a_{3m}=101001000$. By the definition of $X$, to prove $a_1\ldots a_{3m+3}\in \mathcal B^*(X)$ it suffices to prove $a_{3m+1}a_{3m+2}a_{3m+3}\in\set{101, 100}$. Since $11$ is forbidden in $a_1\ldots a_{3\ell}$ by (\ref{eq:forbidden-blocks}), we only need to prove $a_{3m+1}=1$.

 Suppose on the contrary that $a_{3m+1}=0$. By (\ref{eq:sum-one-digits}) and (\ref{eq:forbidden-blocks}) it follows that
 \[ \left(\begin{array}{c}
   d_{3m-8}^1\ldots d_{3m+3}^1\\
   d_{3m-8}^2\ldots d_{3m+3}^2\\
   \overline{d_{3m-8}^\+}\ldots\overline{d_{3m+3}^\+}
 \end{array}\right)=\left(
 \begin{array}{cccc}
   101&001&000&010\\
    010&010&010&100\\
    000&100&101&001
 \end{array}\right).\]
{If $\ell=m+1$, then $k\in\set{3m+3, 3m+4, 3m+5}$, and we consider them separately.
\begin{itemize}
  \item If $k=3m+3$, then by using $d_1^1d_2^1d_3^1=a_1a_2a_3=101$ it follows that $d_{k-1}^1d_k^1d_1^1d_2^1d_3^1=10101\succ a_1\ldots a_5$, leading to a contradiction with (\ref{eq:admissible-condition}).
  \item If $k=3m+4$, then by (\ref{eq:admissible-condition}) we must have $d_k^1=0$. Furthermore, by (\ref{eq:forbidden-blocks}) we must have $\overline{d_k^\+}=0$, which implies $d_k^2=1$ by (\ref{eq:sum-one-digits}). So, $d_k^2d_1^2\ldots d_8^2=101001001\succ a_1\ldots a_9$, leading to a contradiction with (\ref{eq:admissible-condition}).

      \item If $k=3m+5$, then by (\ref{eq:forbidden-blocks}) it follows that $d_{k-1}^1d_k^1=00$. Furthermore, by using $\overline{d_{3m-2}^\+\ldots d_{3m+3}^\+}=101001$ and (\ref{eq:forbidden-blocks}) we must have $\overline{d_{k-1}^\+d_k^\+}=\overline{d_{3m+4}^\+d_{3m+5}^\+}=00$, which implies $d_{k-1}^2d_k^2=11$ by (\ref{eq:sum-one-digits}). This again leads to a contradiction with (\ref{eq:forbidden-blocks}).
\end{itemize}

In the following we assume $\ell>m+1$.    For the next digit $d_{3m+4}^1=a_{3m+4}$ we split the proof into   two subcases.}
\begin{itemize}
  \item[(Ia)]   $d_{3m+4}^1=0$. Note that $\overline{d_{3m-2}^\+\ldots d_{3m+3}^\+}=101001.$   Then by {(\ref{eq:forbidden-blocks})} it follows that $\overline{d_{3m+4}^\+d_{3m+5}^\+d_{3m+6}^\+}=000$. Since $d_{3m+4}^1=0$, by (\ref{eq:sum-one-digits}) we have $d_{3m+4}^2=1$, and thus $d_{3m-1}^2\ldots d_{3m+4}^2=101001$. By {(\ref{eq:forbidden-blocks})} it follows that $d_{3m+5}^2d_{3m+6}^2=00$. Therefore, by (\ref{eq:sum-one-digits}) we obtain that
      \[ \left(\begin{array}{c}
   d_{3m-8}^1\ldots d_{3m+6}^1\\
   d_{3m-8}^2\ldots d_{3m+6}^2\\
   \overline{d_{3m-8}^\+ \ldots d_{3m+6}^\+}
 \end{array}\right)=\left(
 \begin{array}{ccccc}
    101&001&000&010&011\\
    010&010&010&100&100\\
    000&100&101&001&000
 \end{array}\right).\]
 But then $d_{3m+5}^1d_{3m+6}^1=11$,   leading to a contradiction with (\ref{eq:forbidden-blocks}).

 \item[(Ib)]  $d_{3m+4}^1=1$. Then $d_{3m+2}^1d_{3m+3}^1d_{3m+4}^1=101$, and thus by  {(\ref{eq:forbidden-blocks})} it follows that $d_{3m+5}^1d_{3m+6}^1=00$. {By (\ref{eq:sum-one-digits})} this implies that $\overline{d_{3m+4}^\+d_{3m+5}^\+d_{3m+6}^\+}\in\set{001, 010}$, {which yields that
     \[\overline{d_{3m-2}^\+\ldots d_{3m+5}^\+}=10100101\quad\textrm{or}\quad\overline{d_{3m-2}^\+\ldots d_{3m+6}^\+}=101001001,\]
      again leading to a contradiction with (\ref{eq:forbidden-blocks})}.
\end{itemize}
Therefore,
 $a_{3m+1}=1$, and thus $a_{3m+1}a_{3m+2}a_{3m+3}\in\set{101, 100}$.

 Next we prove $d_{3m+1}^2d_{3m+2}^2d_{3m+3}^2=010$. Note that $d_{3m-8}^2\ldots d_{3m}^2=(010)^3$. In view of the choice of $a_{3m+1}a_{3m+2}a_{3m+3}$, if   $a_{3m+1}a_{3m+2}a_{3m+3}=101$, then by (\ref{eq:sum-one-digits}) and (\ref{eq:forbidden-blocks}) we must have
  \[ \left(\begin{array}{c}
   d_{3m-8}^1\ldots d_{3m+3}^1\\
   d_{3m-8}^2\ldots d_{3m+3}^2\\
   \overline{d_{3m-8}^\+ \ldots d_{3m+3}^\+}
 \end{array}\right)=\left(
 \begin{array}{cccc}
   101&001&000&101\\
    010&010&010&010\\
    000&100&101&000
 \end{array}\right),\]
 which gives $d_{3m+1}^2d_{3m+2}^2d_{3m+3}^2=010$. Similarly, if $a_{3m+1}a_{3m+2}a_{3m+3}=100$, then we can also {deduce} that $d_{3m+1}^2d_{3m+2}^2d_{3m+3}^2=010$.

 Case II. $a_{3m-8}\ldots a_{3m}=101000101$. By the definition of $X$, to prove $a_1\ldots a_{3m+3}\in\mathcal B^*(X)$ it suffices to prove $a_{3m+1}a_{3m+2}a_{3m+3}\in\set{000, 001}$. Note that  $d_{3m-2}^1d_{3m-1}^1d_{3m}^1=101$. Then by {(\ref{eq:forbidden-blocks})} it follows that $d_{3m+1}^1d_{3m+2}^1=00$, and thus $a_{3m+1}a_{3m+2}a_{3m+3}\in\set{000, 001}$ as required.

 Next we prove $d_{3m+1}^2d_{3m+2}^2d_{3m+3}^2=010$. If $a_{3m+1}a_{3m+2}a_{3m+3}=000$, then by (\ref{eq:sum-one-digits}) and (\ref{eq:forbidden-blocks}) it follows that
   \[ \left(\begin{array}{c}
   d_{3m-8}^1\ldots d_{3m+3}^1\\
   d_{3m-8}^2\ldots d_{3m+3}^2\\
   \overline{d_{3m-8}^\+ \ldots d_{3m+3}^\+}
 \end{array}\right)=\left(
 \begin{array}{cccc}
   101&000&101&000\\
    010&010&010&010\\
    000&101&000&101
 \end{array}\right),\]
 which gives $d_{3m+1}^2d_{3m+2}^2d_{3m+3}^2=010$. If $a_{3m+1}a_{3m+2}a_{3m+3}=001$, then  $d_{3m+1}^2d_{3m+2}^2d_{3m+3}^2\in\set{010, 100}$. Suppose on the contrary that $d_{3m+1}^2d_{3m+2}^2d_{3m+3}^2=100$. {If $\ell=m+1$, then $k\in\set{3m+3, 3m+4, 3m+5}$, and consider them separately.
 \begin{itemize}
   \item If $k=3m+3$, then \[ \left(\begin{array}{c}
   d_{3m-8}^1\ldots d_{3m+3}^1\\
   d_{3m-8}^2\ldots d_{3m+3}^2\\
   \overline{d_{3m-8}^\+ \ldots d_{3m+3}^\+}
 \end{array}\right)=\left(
 \begin{array}{cccc}
    101&000&101&001 \\
    010&010&010&100 \\
    000&101&000&010
 \end{array}\right).\]
 This implies that $d_k^1d_1^1=11\succ a_1a_2$, leading to a contradiction with (\ref{eq:admissible-condition}).

 \item If $k=3m+4$, then $d_k^1=0$ by (\ref{eq:forbidden-blocks}), and thus $d_{k-1}^1d_k^1d_1^1d_2^1d_3^1=10101\succ a_1\ldots a_5$, also leading to a contradiction with (\ref{eq:admissible-condition}).

 \item If $k=3m+5$, then $d_{k-1}^1d_k^1=00$ by (\ref{eq:forbidden-blocks}), and hence, $d_{k-7}^1\ldots d_k^1d_1^1=101001001\succ a_1\ldots a_9$, leading to a contradiction with (\ref{eq:admissible-condition}).
 \end{itemize}

 In the following we assume $\ell>m+1$.}
 Then by (\ref{eq:sum-one-digits}) and (\ref{eq:forbidden-blocks}) it follows that
       \[ \left(\begin{array}{c}
   d_{3m-8}^1\ldots d_{3m+6}^1\\
   d_{3m-8}^2\ldots d_{3m+6}^2\\
   \overline{d_{3m-8}^\+ \ldots d_{3m+6}^\+}
 \end{array}\right)=\left(
 \begin{array}{ccccc}
    101&000&101&001&000\\
    010&010&010&100&010\\
    000&101&000&010&101
 \end{array}\right)\]
 or
 \[
 \left(\begin{array}{c}
   d_{3m-8}^1\ldots d_{3m+6}^1\\
   d_{3m-8}^2\ldots d_{3m+6}^2\\
   \overline{d_{3m-8}^\+ \ldots d_{3m+6}^\+}
 \end{array}\right)=\left(
 \begin{array}{ccccc}
    101&000&101&001&000\\
    010&010&010&100&101\\
    000&101&000&010&010
 \end{array}\right).\]
 In the first case, we have $\overline{d_{3m+2}^\+\ldots d_{3m+6}^\+}=10101$, and in the second case we also have $d_{3m-1}^2\ldots d_{3m+6}^2=10100101$. Both cases lead to a contradiction {with (\ref{eq:forbidden-blocks})}. Therefore, $d_{3m+1}^2d_{3m+2}^2d_{3m+3}^2=010$.

 Case III. $a_{3m-8}\ldots a_{3m}=101000100$.   To prove $a_1\ldots a_{3m+3}\in \mathcal B^*(X)$ it suffices to prove $a_{3m+1}a_{3m+2}a_{3m+3}=101$. First we prove $a_{3m+1}=1$. If $a_{3m+1}=0$, then by (\ref{eq:sum-one-digits}) and (\ref{eq:forbidden-blocks}) we have
  \[
  \left(\begin{array}{l}
  d_{3m-8}^1\ldots d_{3m+1}^1\\
   d_{3m-8}^2\ldots d_{3m+1}^2\\
   \overline{d_{3m-8}^\+ \ldots  d_{3m+1}^\+}
   \end{array}\right)
   =\left(
   \begin{array}{cccccc}
   101&000&100&0\\
   010&010&010&1\\
   000&101&001&0\end{array}\right).
  \]
  Note that $d_{3m-1}^2d_{3m}^2d_{3m+1}^2=101$ and $\overline{d_{3m-5}^\+\ldots d_{3m+1}^\+}=1010010$. By {(\ref{eq:forbidden-blocks})}  it follows that
  \[
  d_{3m+2}^2d_{3m+3}^2=00=\overline{d_{3m+2}^\+d_{3m+3}^\+}.
  \]
  This together with (\ref{eq:sum-one-digits}) implies $d_{3m+2}^1d^1_{3m+3}=11$, leading to a contradiction with (\ref{eq:forbidden-blocks}). So, $a_{3m+1}=1$. Therefore, by (\ref{eq:sum-one-digits}) and (\ref{eq:forbidden-blocks}) we can deduce that
  \[
  \left(\begin{array}{l}
  d_{3m-8}^1\ldots d_{3m+3}^1\\
   d_{3m-8}^2\ldots d_{3m+3}^2\\
   \overline{d_{3m-8}^\+ \ldots  d_{3m+3}^\+}
   \end{array}\right)
   =\left(
   \begin{array}{cccccc}
   101&000&100&101\\
   010&010&010&010\\
   000&101&001&000\end{array}\right).
  \]
This proves $a_{3m+1}a_{3m+2}a_{3m+3}=101$ and $d_{3m+1}^2d_{3m+2}d_{3m+3}^2=010$ as desired.

Case IV. $a_{3 m-8}\ldots a_{3m}=100101001$. To prove $a_1\ldots a_{3m+3}\in\mathcal B^*(X)$ it suffices to prove $a_{3m+1}a_{3m+2}a_{3m+3}=000$. Note that $a_{3m-5}\ldots a_{3m}=101001$ {and $a_1\ldots a_{3m+3}=d_1^1\ldots d_{3m+3}^1$}. By {(\ref{eq:forbidden-blocks})} it follows that $a_{3m+1}a_{3m+2}a_{3m+3}=000$. Furthermore, by (\ref{eq:sum-one-digits}) and (\ref{eq:forbidden-blocks}) we must have
\[
  \left(\begin{array}{l}
  d_{3m-8}^1\ldots d_{3m+3}^1\\
   d_{3m-8}^2\ldots d_{3m+3}^2\\
   \overline{d_{3m-8}^\+ \ldots  d_{3m+3}^\+}
   \end{array}\right)
   =\left(
   \begin{array}{cccccc}
   100&101&001&000\\
   010&010&010&010\\
   001&000&100&101\end{array}\right),
  \]
  which implies $d_{3m+1}^2d_{3m+2}^2d_{3m+3}^2=010$.

Case V. $a_{3m-8}\ldots a_{3m}=100101000$. To prove $a_1\ldots a_{3m+3}\in\mathcal B^*(X)$ it suffices to prove $a_{3m+1}a_{3m+2}a_{3m+3}\in\set{100, 101}$. In view of (\ref{eq:forbidden-blocks}) and $a_1\ldots a_{3m+3}=d_1^1\ldots d_{3m+3}^1$ it suffices to prove $a_{3m+1}=1$.

Suppose on the contrary that $a_{3m+1}=0$. Then by (\ref{eq:sum-one-digits}) and (\ref{eq:forbidden-blocks}) it follows that
\[
  \left(\begin{array}{l}
  d_{3m-8}^1\ldots d_{3m+3}^1\\
   d_{3m-8}^2\ldots d_{3m+3}^2\\
   \overline{d_{3m-8}^\+ \ldots d_{3m+3}^\+}
   \end{array}\right)
   =\left(
   \begin{array}{cccccc}
   100&101&000&010\\
   010&010&010&100\\
   001&000&101&001\end{array}\right).
  \]
  {If $\ell=m+1$, then $k\in\set{3m+3, 3m+4, 3m+5}$. By the same argument as in Case I we can always deduce a contradiction.}
%
%
So in the following we assume $\ell>m+1$, and we consider two subcases.
  \begin{itemize}
  \item[{\rm(Va)}] $a_{3m+4}=1$. Then by (\ref{eq:sum-one-digits}) and (\ref{eq:forbidden-blocks}) we must have
  \[
  \left(\begin{array}{l}
  d_{3m-8}^1\ldots d_{3m+6}^1\\
   d_{3m-8}^2\ldots d_{3m+6}^2\\
   \overline{d_{3m-8}^\+ \ldots  d_{3m+6}^\+}
   \end{array}\right)
   =\left(
   \begin{array}{ccccccc}
   100&101&000&010&101\\
   010&010&010&100&010\\
   001&000&101&001&000\end{array}\right),
  \]
which implies that $d^1_{3m+2}\ldots d_{3m+6}^1=10101$, leading to a contradiction {with (\ref{eq:forbidden-blocks})}.

\item[{\rm(Vb)}] $a_{3m+4}=0$. Then by (\ref{eq:sum-one-digits}) and (\ref{eq:forbidden-blocks}) it {follows} that
  \[
  \left(\begin{array}{l}
  d_{3m-8}^1\ldots d_{3m+6}^1\\
   d_{3m-8}^2\ldots d_{3m+6}^2\\
   \overline{d_{3m-8}^\+}\ldots \overline{d_{3m+6}^\+}
   \end{array}\right)
   =\left(
   \begin{array}{ccccccc}
   100&101&000&010&010\\
   010&010&010&100&101\\
   001&000&101&001&000\end{array}\right).
  \]
  This yields that $d_{3m-1}^2\ldots d_{3m+6}^2=10100101$, again leading to a contradiction {with (\ref{eq:forbidden-blocks})}.
\end{itemize}
Therefore, $a_{3m+1}=1$, and then $a_{3m+1}a_{3m+2}a_{3m+3}\in\set{101, 100}$. If $a_{3m+1}a_{3m+2}a_{3m+3}=101$, then by (\ref{eq:sum-one-digits}) and (\ref{eq:forbidden-blocks}) we must have
\[
  \left(\begin{array}{l}
  d_{3m-8}^1\ldots d_{3m+3}^1\\
   d_{3m-8}^2\ldots d_{3m+3}^2\\
   \overline{d_{3m-8}^\+ \ldots  d_{3m+3}^\+}
   \end{array}\right)
   =\left(
   \begin{array}{cccccc}
   100&101&000&101\\
   010&010&010&010\\
   001&000&101&000\end{array}\right).
  \]
This proves $d_{3m+1}^2d_{3m+2}^2d_{3m+3}^2=010$. {Similarly}, for $a_{3m+1}a_{3m+2}a_{3m+3}={100}$ we can also deduce that $d_{3m+1}^2d_{3m+2}^2d_{3m+3}^2=010$.

\medskip

Hence, by Cases I--V  it follows that $a_1\ldots a_{3m+3}\in\mathcal B^*(X)$ and $d_{1}^2\ldots d_{3m+3}^2=(010)^{m+1}$.
Proceeding this argument we can prove that  $a_1\ldots a_{3\ell}\in\mathcal B^*(X)$ and $d_1^2\ldots d_{3\ell}^2=(010)^\ell$, completing  the proof.
\end{proof}

 Take $k=3(2m+1)\cdot 2^n$ with $m\in\N_{\ge 2}$ and $n\in\N_0$. Now we start to prove the lower bound of $\beta_k$. Write
 \[\bm\la:=\la_1\ldots\la_{3\cdot 2^n}\quad\textrm{and}\quad  \bm\ga:=\ga_1\ldots \ga_{3\cdot 2^n}=\Theta(\bm\la).\]
 Then
 \[
 (\la_1\ldots \la_{3(2^{n+1}+2^n)}(\la_1\ldots\la_{3\cdot 2^{n+1}}^-)^{m-1})^\f=(\bm\la\bm\ga^+\bm\ga(\bm\la\bm\ga)^{m-1})^\f,\quad\de(\beta_c)=\bm\la\bm\ga^+\bm\ga\bm\la\bm\ga\bm\la^-\bm\la\ldots.
 \]
In view of (\ref{eq:etak-lower-bound}),    it suffices to prove that there is no admissible block $a_1\ldots a_k=:\a_1\ldots \a_{2m+1}$ satisfying
 \begin{equation}\label{eq:etak-lower-bound-check}
 \bm\la\bm\ga^+\bm\ga\bm\la\bm\ga\bm\la^-\bm\la\ldots=\de(\beta_c)\lle (\a_1\ldots \a_{2m+1})^\f\prec (\bm{\la\ga^+\ga(\la\ga)}^{m-1})^\f,
 \end{equation}where each $\a_i$ is a block in $\set{0,1}^*$ of length $3\cdot 2^n$. Then {$\bm\la, \bm\ga$ and} each block $\a_i$ have the same  length. For a {block} $\d=d_1\ldots d_s\in  \Omega^s$ we write $\d^1:=d_1^1\ldots d_s^1, \d^2:=d_1^2\ldots d_s^2$ and $\overline{\d^\+}:=\overline{d_1^\+}\ldots \overline{d_s^\+}$.

  Some ideas on proving  the lower bound of $\beta_k$ are illustrated in the following lemma.

 \begin{lemma}
   \label{lem:lower-bound-234}
   Let $k=3(2m+1)2^n$ with $m\in\set{2,3,4}$ and $n\in\N_0$. Then
   \[
   \de(\beta_k)=(\bm{\la\ga^+\ga(\la\ga)}^{m-1})^\f.
   \]
 \end{lemma}
 \begin{proof}
   By Lemma \ref{lem:upper-bound-eta-3k} we only need to prove that there is no admissible {block} $\a_1\ldots\a_{2m+1}$ satisfying (\ref{eq:etak-lower-bound-check}).
   It is easy to check this for $m=2$.

   For $m=3$ we assume on the contrary that there is an admissible  {block} $\a_1\ldots \a_7$ satisfying  (\ref{eq:etak-lower-bound-check}), i.e., $
   \bm{\la\ga^+\ga\la\ga\la^-\la}\lle \a_1\ldots \a_7\prec \bm{\la\ga^+\ga\la\ga\la\ga}.$ This
    implies that
    \[\a_1\ldots \a_5=\bm{\la\ga^+\ga\la\ga} \quad\textrm{and}\quad \a_6\in\set{\bm\la^-, \bm\la}.\]
     Since $\a_1\ldots \a_7$ is admissible and it has a prefix $\la_1\ldots\la_{9}=101001000$, by Proposition \ref{prop:key-prop} it follows that there exists a unique {block} $d_1\ldots d_{3(2m+1)2^n}=\d_1\ldots\d_{2m+1}$ with each $\d_i\in \Omega^{3\cdot 2^n}$
   such that
   \[
   \left(
   \begin{array}
     {c}
     \d_1^1\ldots\d_7^1\\
     \d_1^2\ldots \d_7^2\\
     \overline{\d_1^\+}\ldots \overline{\d_7^\+}
   \end{array}\right)=
   \left(
   \begin{array}
     {ccccccc}
     \bm\la&\bm\ga^+&\bm\ga&\bm\la&\bm\ga&\a_6&\a_7\\
     \bm\nu&\bm\nu&\bm\nu&\bm\nu&\bm\nu&\bm\nu&\bm\nu\\
     \bm\ga&\bm\la^-&\bm\la&\bm\ga&\bm\la&\Theta(\a_6)&\Theta(\a_7)
   \end{array}\right),
   \]
where we reserve the notation $\bm\nu:=(010)^{2^n}$.
    If $\a_6=\bm\la^-$, then by using the admissibility of $\a_1\ldots \a_7$ it follows that
    {\[
    \a_7\bm\la=\a_7\a_1\lle\a_1\la_2=\bm{\la\ga^+},
    \]
    which together with $\bm\la\succ\bm\ga^+$ implies} $\a_7\prec \bm\la$. Thus, $(\a_1\ldots \a_7)^\f\prec \bm{\la\ga^+\ga\la\ga\la^-\la}0^\f\lle \de(\beta_c)$, leading to a contradiction with (\ref{eq:etak-lower-bound-check}). So $\a_6=\bm\la$. By (\ref{eq:etak-lower-bound-check}) we  have $\a_7\prec\bm\ga$. Then by using $\d_7^2=\bm\nu$ and (\ref{eq:sum-one-digits}) it follows that $\overline{\d_7^\+}=\Theta(\a_7)\succ\bm\la=\a_1$, this again leads to a contradiction with (\ref{eq:admissible-condition}).

   Next we consider $m=4$. Suppose on the contrary {that} there is an admissible  {block} $\a_1\ldots \a_9$ satisfying  (\ref{eq:etak-lower-bound-check}). Then
   \begin{equation}\label{eq:m4-1}
   \bm{\la\ga^+\ga\la\ga\la^-\la\ga^+\ga}\lle \a_1\ldots \a_9\prec \bm{\la\ga^+\ga\la\ga\la\ga\la\ga}.
   \end{equation}
   So, $\a_1\ldots \a_5=\bm{\la\ga^+\ga\la\ga}$. By Proposition \ref{prop:key-prop} there exists a unique block $\d_1\ldots \d_9$ with each $\d_i\in \Omega^{3\cdot 2^n}$ such that
    \[
   \left(
   \begin{array}
     {c}
     \d_1^1\ldots\d_9^1\\
     \d_1^2\ldots \d_9^2\\
     \overline{\d_1^\+ \ldots  \d_9^\+}
   \end{array}\right)=
   \left(
   \begin{array}
     {ccccccccc}
     \bm\la&\bm\ga^+&\bm\ga&\bm\la&\bm\ga&\a_6&\a_7&\a_8&\a_9\\
     \bm\nu&\bm\nu&\bm\nu&\bm\nu&\bm\nu&\bm\nu&\bm\nu&\bm\nu&\bm\nu\\
     \bm\ga&\bm\la^-&\bm\la&\bm\ga&\bm\la&\Theta(\a_6)&\Theta(\a_7)&\Theta(\a_8)&\Theta(\a_9)
   \end{array}\right).
   \]
   Note by (\ref{eq:m4-1}) that $\a_6\in\set{\bm\la^-, \bm\la}$. If $\a_6=\bm\la^-$, then by (\ref{eq:admissible-condition}) we have
   \[\a_7\a_8\a_9\bm{\la\ga^+}=\a_7\a_8\a_9\a_1\a_2\lle\a_1\ldots \a_5=\bm{\la\ga^+\ga\la\ga},\] which implies $\a_7\a_8\a_9\prec\a_1\a_2\a_3= \bm{\la\ga^+\ga}$, and thus $(\a_1\ldots \a_9)^\f\prec \de(\beta_c)$, leading to a contradiction {with (\ref{eq:m4-1})}. Thus, $\a_1\ldots \a_6=\bm{\la\ga^+\ga\la\ga\la}$. By (\ref{eq:m4-1}) we have $\a_7\lle\bm\ga$. If $\a_7\prec \bm\ga$, then $\overline{\d_7^\+}=\Theta(\a_7)\succ \bm\la= \a_1$, leading to a contradiction with the admissibility of $\a_1\ldots \a_9$. So, $\a_7=\bm\ga$. Again by (\ref{eq:m4-1}) we have $\a_8\lle \bm\la$.

   If $\a_8\prec \bm\la$, then we can deduce that $\overline{\d_8^\+}=\Theta(\a_8)\lge\bm\ga^+$, and thus $\overline{\d_7^\+ \d_8^\+}=\Theta(\a_7)\Theta(\a_8)\lge\bm{\la\ga^+}=\a_1\a_2$. Note by (\ref{eq:admissible-condition}) that $\a_9\bm{\la}=\a_9\a_1\lle  {\a_1\a_2}=\bm{\la\ga^+}$. Then by using $\bm\la\succ \bm\ga^+$ we must have $\a_9\prec \bm\la$.  Thus,  we can deduce that \[\overline{\d_7^\+\d_8^\+\d_9^\+}=\Theta(\a_7 \a_8 \a_9)\succ\bm{\la\ga^+\ga}=\a_1\a_2\a_3,\] leading to a contradiction with (\ref{eq:admissible-condition}). Therefore, we must have $\a_8=\bm\la$. This together with (\ref{eq:m4-1}) implies that $\a_9\prec \bm\ga$, and then $\overline{\d_9^\+}=\Theta(\a_9)\succ \bm\la=\a_1$, a contradiction.
 \end{proof}

\begin{proof}
  [Proof of Theorem \ref{th:eta-period-3n}]
  By Proposition \ref{prop:m=1},  Lemmas \ref{lem:upper-bound-eta-3k} and \ref{lem:lower-bound-234} it suffices to prove for $m\in\N_{\ge 5}$ that   no admissible {block} $\a_1\ldots \a_{2m+1}$   satisfies (\ref{eq:etak-lower-bound-check}).

  Suppose on the contrary {that} there is an admissible {block} $\a_1\ldots\a_{2m+1}$ satisfying (\ref{eq:etak-lower-bound-check}), i.e.,
   \begin{equation*}
 \bm{\la\ga^+\ga\la\ga\la^-\la\ga^+}\ldots= \de(\beta_c)\lle (\a_1\ldots \a_{2m+1})^\f\prec (\bm{\la\ga^+\ga(\la\ga)}^{m-1})^\f.
  \end{equation*}
 Then $\a_1\ldots \a_5=\bm{\la\ga^+\ga\la\ga}$ and $\a_6\in\set{\bm\la^-, \bm\la}$. First we prove $\a_6=\la$.

 If $\a_6=\bm\la^-$, then by Proposition \ref{prop:key-prop} there exists a unique block $\d_1\ldots \d_{2m+1}$ with each $\d_i\in \Omega^{3\cdot 2^n}$ such that
   \[
   \left(
   \begin{array}
     {c}
     \d_1^1\ldots\d_{2m+1}^1\\
     \d_1^2\ldots \d_{2m+1}^2\\
     \overline{\d_1^\+}\ldots \overline{\d_{2m+1}^\+}
   \end{array}\right)=
   \left(
   \begin{array}
     {ccccccccc}
     \bm\la&\bm\ga^+&\bm\ga&\bm\la&\bm\ga&\bm\la^-&\a_7&\ldots&\a_{2m+1}\\
     \bm\nu&\bm\nu&\bm\nu&\bm\nu&\bm\nu&\bm\nu&\bm\nu&\ldots&\bm\nu\\
     \bm\ga&\bm\la^-&\bm\la&\bm\ga&\bm\la&\bm\ga^+&\Theta(\a_7)&\ldots&\Theta(\a_{2m+1})
   \end{array}\right).
   \]
   Note that $\Theta(\a_i\a_{i+1})=\overline{\d_i^\+\d_{i+1}^\+}\lle \a_1\a_2=\bm{\la\ga^+}$ for all $1\le i\le 2m$. Then
    \begin{equation}
      \label{eq:lower-m-2}
      \a_i\a_{i+1}\lge \Theta(\bm\la)\Theta(\bm\ga^+)=\bm{\ga\la^-}.
    \end{equation}
    Furthermore, by using the admissibility of $\a_1\ldots\a_{2m+1}$ we have
   \begin{equation}
     \label{eq:lower-m-3}
     \a_{i+1}\ldots \a_{i+6}\lle \a_1\ldots \a_6=\bm{\la\ga^+\ga\la\ga\la^-}\quad\forall 1\le i\le 2m-5.
   \end{equation}
   Note also that $\bm{\la\ga^+}\Theta(\a_7\ldots \a_{10})=\overline{\d_5^\+\ldots \d_{10}^\+}\lle \a_1\ldots \a_6=\bm{\la\ga^+\ga\la\ga\la^-}$. Then
   \begin{equation}
     \label{eq:lower-m-4}
     \a_7\ldots \a_{10}\lge \bm{\la\ga\la\ga^+}.
   \end{equation}
   By (\ref{eq:lower-m-3}) and (\ref{eq:lower-m-4}) it follows that $\a_7\a_8\in\set{\bm{\la\ga}, \bm{\la\ga^+}}$. But if we take $\a_7\a_8=\bm{\la\ga}$, then  (\ref{eq:lower-m-3}) and (\ref{eq:lower-m-4}) imply that $\a_9\a_{10}=\bm{\la\ga^+}$. Suppose $\a_7\a_8=\bm{\la\ga^+}$. Then by (\ref{eq:lower-m-2}) and (\ref{eq:lower-m-3}) it follows that $\a_9\a_{10}\in\set{\bm{\ga\la^-}, \bm{\ga\la}}$. Once we take $\a_9\a_{10}=\bm{\ga\la}$, then we must have $\a_{11}\a_{12}=\bm{\ga\la^-}$ by (\ref{eq:lower-m-2}) and (\ref{eq:lower-m-3}).

   Proceeding this argument we can conclude that the {block} $\a_1\ldots \a_{2m+1}$ can be represented via the directed graph $X_{\mathcal G}$ as plotted in Figure \ref{fig:2}.
    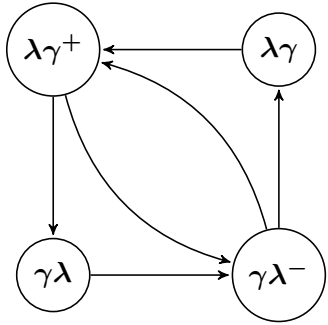
\begin{figure}[h!]
  \centering
  \begin{tikzpicture}[->,>=stealth',shorten >=1pt,auto,node distance=3cm, semithick,scale=3]

  \tikzstyle{every state}=[minimum size=0pt,fill=none,draw=black,text=black]

  \node[state] (A)                    {$\bm{\ga\la}$ };
  \node[state]         (B) [ right of=A] { $\bm{\ga\la^-}$};
  \node[state]         (C) [ above of=A] {$\bm{\la\ga^+}$};
  \node[state](D)[right of=C]{$\bm{\la\ga}$};

  \path[->,every loop/.style={min distance=0mm, looseness=25}]
  (C) edge[->,left] node{} (A)
  (C) edge[bend right,->,right] node{} (B)
(D) edge[->,above] node{} (C)
(A)edge[->,below] node{} (B)
(B) edge[bend right,->,right] node{} (C)
(B) edge[->,right] node{} (D)
;
\end{tikzpicture}
\caption{The directed graph representing the possible connections of blocks $\bm{\la\ga^+}, \bm{\la\ga}, \bm{\ga\la}$ and $\bm{\ga\la^+}$.}
\label{fig:2}
\end{figure}
Therefore, $\a_1\ldots \a_{2m+1}$ must  {have a suffix} in
\[
\set{\bm{\la\ga^+\ga},\quad \bm{\ga\la\ga}, \quad \bm{\ga\la^-\la},\quad \bm{\la\ga\la}}.
\]
If $\a_{2m-1}\a_{2m}\a_{2m+1}=\bm{\la\ga^+\ga}$, then $\a_{2m-1}\a_{2m}\a_{2m+1}\a_1\a_2=\bm{\la\ga^+\ga\la\ga^+}\succ\a_1\ldots\a_5$, contradicting to the admissibility of $\a_1\ldots\a_{2m+1}$. If $\a_{2m-1}\ldots\a_{2m+1}=\bm{\ga\la\ga}$, then we must have
\[\a_{2m-3}\ldots \a_{2m+1}\a_1=\bm{\la\ga^+\ga\la\ga\la}\succ\bm{\la\ga^+\ga\la\ga\la^-}=\a_1\ldots \a_6,\] again leading to a contradiction. If $\a_{2m-1}\a_{2m}\a_{2m+1}\in\set{\bm{\ga\la^-\la}, \bm{\la\ga\la}}$, then
$\a_{2m+1}\a_1=\bm{\la\la}\succ\bm{\la\ga^+}=\a_1\a_2$, a contradiction.

Therefore, $\a_6=\bm\la$. By (\ref{eq:etak-lower-bound-check}) and (\ref{eq:lower-m-2}) it follows  that $\a_7\a_8\in\set{\bm{\ga\la^-}, \bm{\ga\la}}$. If $\a_7\a_8=\bm{\ga\la^-}$, then by the same argument as above for $\a_5\a_6=\bm{\ga\la^-}$ it will lead to a contradiction. This means $\a_7\a_8=\bm{\ga\la}$. Proceeding this argument and using (\ref{eq:etak-lower-bound-check}),  (\ref{eq:lower-m-2}) we will either deduce a contradiction  if $\a_{2i-1}\a_{2i}=\bm{\ga\la^-}$ for some $i$, or will eventually have
\[
\a_1\ldots\a_{2m+1}=\bm{\la\ga^+\ga(\la\ga)}^{m-1},
\]
which again leads to a contradiction with (\ref{eq:etak-lower-bound-check}). This completes the proof.
\end{proof}

 Let $(\tau_i)_{i=0}^\f=01101001\ldots$ be the classical Thue-Morse sequence {(cf.~\cite{Allouche-Shallit-1999-1})}. Note that our definition of $(\la_i)_{i=1}^\f$ in (\ref{eq:la-ga}) is related to the sequence $(\tau_i)_{i=0}^\f$:
\begin{equation}\label{eq:relation-tau-lambda}
\la_{3n+1}=\tau_{2n+1},\quad\la_{3n+2}=0,\quad \la_{3n+3}=\tau_{2n+2}\quad\forall n\ge 0.
\end{equation}
Recall the Sharkovskii ordering defined in (\ref{eq:sharkovski}).
To complete our proof of Theorem \ref{th:main result} (ii) we recall the following result from \cite{Allouche-Clarke-Sidorov-2009}.
\begin{proposition}
  \label{prop:Sharkovski-order}
Let $k=(2m+1)2^n$ with $m, n\in\N_0$, and let
\[
\xi_k=\left\{
\begin{array}
  {lll}
  (\tau_1\ldots\tau_{2^n}^-)^\f&\textrm{if}& m=0,\\
  (\tau_1\ldots \tau_{2^{n+1}+2^n}(\tau_1\ldots \tau_{2^{n+1}}^-)^{m-1})^\f&\textrm{if}& m\ge 1.
\end{array}\right.
\]
Then $\xi_k\succ\xi_{\ell}$ if and only if $k\rhd \ell$ in {the} Sharkovskii ordering.
\end{proposition}

\begin{proof}
  [Proof of Theorem \ref{th:main result} (i)--(ii)]It is clear that $\beta_1=1$. Then
  (i) follows easily from Theorem \ref{th:lower-higher-orders} (iii). For (ii),
  the formula for $\de(\beta_{k})$ with $k\in 3\N$ follows directly from Theorem \ref{th:lower-higher-orders} (i) and Theorem \ref{th:eta-period-3n}. Here we point out that we have used the convention  $\de(\beta_1)=0^\f$. {For the Sharkovskii ordering statement, by (\ref{eq:relation-tau-lambda}),   Theorem \ref{th:lower-higher-orders} (i) and Theorem \ref{th:eta-period-3n} it follows that
  \[
  \de(\beta_{3k})\succ\de(\beta_{3\ell})\quad\Longleftrightarrow\quad \xi_{k}\succ\xi_{\ell}.
  \]
  Then by Lemma \ref{lem:delta-beta} and Proposition \ref{prop:Sharkovski-order} we conclude that $\beta_{3k}>\beta_{3\ell}$ if and only if $k\rhd \ell$ in the Sharkovskii ordering.}
 \end{proof}

 \section{Minimal base for $\U_\beta$ to have a sequence of period  $3\N+2$}\label{sec:3N+2}
 In this section we will consider   $\beta_{3\ell+2}$ for any $\ell\in\N$, and prove Theorem \ref{th:main result} (iv).
 First we prove it   for $\ell=1,2,3$.
\begin{lemma}
  \label{lem:below-3}
 If $\ell\in\set{1,2,3}$, then  $\de(\beta_{3\ell+2})=(101(001)^{\ell-1}00)^\f$.
\end{lemma}
\begin{proof}
First we consider  $\ell=1$. Let $\de(\beta_5)=(a_1\ldots a_5)^\f$. Then by Theorem \ref{th:lower-higher-orders}  it follows that
\[
101001\ldots=\de(\beta_c)\lle(a_1\ldots a_5)^\f\lle \de(\beta_2)=(10)^\f.
\]This implies $a_1\ldots a_4=1010$. {Note that $a_1\ldots a_5$ is admissible.}
By (\ref{eq:admissible-condition}) it follows that $a_1\ldots a_5=10100$ as desired.

Next we consider $\ell=2$. Let $\de(\beta_8)=(a_1\ldots a_8)^\f$.  Then by Theorem \ref{th:lower-higher-orders} we obtain that
\[
10100100\ldots\lle(a_1\ldots a_8)^\f\lle (10)^\f.
\]
This implies that $a_1\ldots a_8\lge 10100100$. On the other hand, note that \[d_1\ldots d_8=\al_1\al_2\al_1\al_0\al_2\al_1\al_2\al_0\in \Omega^8\] is an aperiodic {block} satisfying (\ref{eq:admissible-condition}) and $d_1^1\ldots d_8^1=10100100$. So, {the block $10100100$ is admissible, and thus} $a_1\ldots a_8=10100100$.

Finally, we consider  $\ell=3$.
 Suppose $\de(\beta_{11})=(a_1\ldots a_{11})^\f$. Note that $d_1\ldots d_{11}=\al_1\al_2\al_1(\al_0\al_2\al_1)^2\al_0\al_2\in \Omega^{11}$ is an aperiodic {block} satisfying (\ref{eq:admissible-condition}) and $d_1^1\ldots d_{11}^1=101(001)^2 00$. Then by Theorem \ref{th:lower-higher-orders} it follows that
  \begin{equation}\label{eq:l=3-1}
  101001000101\ldots=\de(\beta_c)\lle (a_1\ldots a_{11})^\f\lle(101(001)^200)^\f.
  \end{equation}
  This implies that $a_1\ldots a_8=10100100$. By (\ref{eq:admissible-condition}) and (\ref{eq:l=3-1}) it suffices to prove $a_9=1$. Suppose on the contrary that $a_9=0$. {Then by (\ref{eq:l=3-1}) and Lemma \ref{lem:delta-beta} it follows that $a_1\ldots a_{11}=10100100010$. This yields that
  \[
 a_{10}a_{11}a_1a_2a_3=10101\succ a_1\ldots a_5,
  \]
  leading to a contradiction with (\ref{eq:admissible-condition}).} This proves $a_9=1$, and  completes the proof.
\end{proof}
In the following we assume $\ell\ge 4$.
 In the following proposition we give the bounds on $\beta_{3\ell+2}$   for     $\ell\in\N$.

\begin{proposition}
  \label{prop:above eta-9}\mbox{}
    For any $\ell\in\N$, we have
    \[(101001000)\prec\de(\beta_{3\ell+2})\lle (101(001)^{\ell-1}00)^\f.\]
    \end{proposition}
\begin{proof}
By Lemma \ref{lem:below-3} it suffices to consider  $\ell\ge 4$.  Take
\[
(d_1\ldots d_{3\ell+2})^\f=(\al_1\al_2\al_1(\al_0\al_2\al_1)^{\ell-1}\al_2\al_0)^\f.
\]
Then
\begin{align*}
  (d_1^1\ldots d_{3\ell+2}^1)^\f&=(101(001)^{\ell-1}00)^\f,\\
  (d_1^2\ldots d_{3\ell+2}^2)^\f&=(010(010)^{\ell-1}10)^\f,\\
  (\overline{d_1^\+}\ldots \overline{d_{3\ell+2}^\+})^\f&=(000(100)^{\ell-1}01)^\f.
\end{align*}
 By Lemma \ref{lem:lexicographical-periodic-U} it follows that $(d_1\ldots d_{3\ell+2})^\f\in\U_\beta$ if $\de(\beta)\succ (101(001)^{\ell-1}00)^\f$. This proves $\de(\beta_{3\ell+2})\lle (101(001)^{\ell-1}00)^\f$.

Next we prove $\de(\beta_{3\ell+2})\succ (101001000)^\f$. Write   $\delta(\beta_{3\ell+2})= (a_1\ldots a_{3\ell+2})^\f$. {Then $a_1\ldots a_{3\ell+2}$ is admissible. So there exists an aperiodic block $d_1\ldots d_{3\ell+2}\in \Omega^{3\ell+2}$ satisfying (\ref{eq:admissible-condition}) and $d_1^1\ldots d_{3\ell+2}^1=a_1\ldots a_{3\ell+2}$.}
  Suppose on the contrary that $\de(\beta_{3\ell+2})\lle (101001000)^\f$. Then by Theorem \ref{th:lower-higher-orders} it follows that
\[
101001000101\ldots=\de(\beta_c)\lle (a_1\ldots a_{3\ell+2})^\f\lle (101001000)^\f,
\] which implies $a_1\ldots a_9=101001000$. {Since $d_1\ldots d_{3\ell+2}$ is a representation of $a_1\ldots a_{3\ell+2}$, by {Proposition \ref{prop:key-prop}}}
  it follows that $a_1\ldots a_{3\ell}\in\mathcal B^*(X)$ and
\begin{equation}
  \label{eq:(3n+2)}
  d_1^1\ldots d_{3\ell}^1=a_1\ldots a_{3\ell},\quad d_1^2\ldots d_{3\ell}^2=(010)^\ell\quad\textrm{and}\quad \overline{d_1^\+\ldots d_{3\ell}^\+}=\Theta(a_1\ldots a_{3\ell}).
\end{equation}
Note by (\ref{eq:admissible-condition}) that $a_{3\ell+1}a_{3\ell+2}=00$.
We will finish the proof by considering the following four cases: $a_{3\ell-2}a_{3\ell-1}a_{3\ell}\in\set{101, 001, 000, 100}$.
\begin{itemize}
  \item If $a_{3\ell-2}a_{3\ell-1}a_{3\ell}=101$, then $d^1_{3\ell-2}\ldots d^1_{3\ell+2}d_1^1d_2^1d_3^1=10100101\succ a_1\ldots a_8$, leading to a contradiction with (\ref{eq:admissible-condition}).
   \item If $a_{3\ell-2}a_{3\ell-1}a_{3\ell}=001$, then by   $a_1\ldots a_{3\ell}\in\mathcal B^*(X)$ it follows that {$a_{3\ell-5}a_{3\ell-4}a_{3\ell-3}=101$, and thus} $d^1_{3\ell-5}\ldots d^1_{3\ell+2}d_1^1=101001001\succ a_1\ldots a_9$, again leading to a contradiction with (\ref{eq:admissible-condition}).

       \item If $a_{3\ell-2}a_{3\ell-1}a_{3\ell}=000$, then $a_{3\ell-2}\ldots a_{3\ell+2}=00000$, and thus by (\ref{eq:(3n+2)}) it follows that
       \[
       \left(
       \begin{array}
         {c}
         d_{3\ell-2}^1\ldots d_{3\ell+2}^1\\
         d_{3\ell-2}^2\ldots d_{3\ell+2}^2\\
         \overline{d_{3\ell-2}^\+ \ldots  d_{3\ell+2}^\+}
       \end{array}\right)=
       \left(
       \begin{array}
         {cc}
         000&00\\
         010&**\\
         101&**
       \end{array}\right).
       \]
       Note that $\overline{d_{3\ell-2}^\+}\ldots \overline{d_{3\ell+2}^\+}\lle a_1\ldots a_5=10100$. Then $\overline{d_{3\ell+1}^\+}\overline{d_{3\ell+2}^\+}=00$, and thus {by (\ref{eq:sum-one-digits})} $d_{3\ell+1}^2d_{3\ell+2}^2=11$, again leading to a contradiction.

       \item If $a_{3\ell-2}a_{3\ell-1}a_{3\ell}=100$, then by using $a_1\ldots a_{3\ell}\in\mathcal B^*(X)$ it follows that  {$a_{3\ell-5}a_{3\ell-4}a_{3\ell-3}=000$.} Thus, by (\ref{eq:(3n+2)}) we obtain
          \[
       \left(
       \begin{array}
         {c}
         d_{3\ell-5}^1\ldots d_{3\ell+2}^1\\
         d_{3\ell-5}^2\ldots d_{3\ell+2}^2\\
         \overline{d_{3\ell-5}^\+ \ldots  d_{3\ell+2}^\+}
       \end{array}\right)=
       \left(
       \begin{array}
         {ccc}
         000&100&00\\
         010&010&**\\
         101&001&**
       \end{array}\right).
       \]
       Note that $\overline{d_{3\ell-5}^\+}\ldots \overline{d_{3\ell+2}^\+}\lle a_1\ldots a_8=10100100$. Then $\overline{d_{3\ell+1}^\+}\overline{d_{3\ell+2}^\+}=00$, and thus  {by (\ref{eq:sum-one-digits})} $d_{3\ell+1}^2d_{3\ell+2}^2=11$,  again leading to a contradiction.
\end{itemize}
Therefore,  $\de(\beta_{3\ell+2})\succ (101001000)^\f$, completing the proof.
\end{proof}

{Note by Theorem \ref{th:main result} (ii) that $\beta_{3\ell}\le \beta_{9}$ for all $\ell\in\N$, and $\de(\beta_9)=(101001000)^\f$. Then by Proposition \ref{prop:above eta-9} and Lemma \ref{lem:delta-beta} it follows that $\beta_{3\ell+2}>\beta_9$ for all $\ell\in\N$.}
Before proving Theorem \ref{th:main result} (iv) we still need the following lemma.
\begin{lemma}
  \label{lem:extended-admissible}
  Let $\ell\in\N_{\ge 4}$ and $1\le m\le \ell-1$. If $a_1\ldots a_{3\ell+2}$ is admissible with
  \begin{equation}\label{eq:extend-1}
  a_1\ldots a_{3m+3}=101(001)^{m-1}000,
  \end{equation}
  then for any representation $d_1\ldots d_{3\ell+2}\in \Omega^{3\ell+2}$ of $a_1\ldots a_{3\ell+2}$ we have
  \[
  d_1^2\ldots d^2_{3m+3}=(010)^{m+1}\quad\textrm{or}\quad \overline{d_1^\+\ldots d_{3m+3}^\+}=(010)^{m+1}.
  \]
\end{lemma}
We point out that for $m=2$ this lemma {is a special case of  {Proposition \ref{prop:key-prop}}.}
\begin{proof}
  Note by (\ref{eq:extend-1}) and (\ref{eq:admissible-condition}) that
  \begin{equation}
    \label{eq:nov-19}
    11\textrm{ and }10101\textrm{ are forbidden in } d_1^2\ldots d^2_{3m+3}\textrm{ and }\overline{d_1^\+\ldots d_{3m+3}^\+}.
  \end{equation}    Then by {(\ref{eq:sum-one-digits}) and} (\ref{eq:extend-1}) it follows that
  \begin{equation}
    \label{eq:extend-2}
    \left(
    \begin{array}
      {l }
      d_{3m+1}^1 d_{3m+2}^1 d_{3m+3}^1\\
       d_{3m+1}^2 d_{3m+2}^2 d_{3m+3}^2\\
        \overline{d_{3m+1}^\+  d_{3m+2}^\+ d_{3m+3}^\+}
    \end{array}\right)=
    \left(
    \begin{array}
      {lll}
      0&0&0\\
      0&1&0\\
      1&0&1
    \end{array}\right)
  \end{equation}
  or
    \begin{equation}
    \label{eq:extend-3}
     \left(
    \begin{array}
      {l }
      d_{3m+1}^1 d_{3m+2}^1 d_{3m+3}^1\\
       d_{3m+1}^2 d_{3m+2}^2 d_{3m+3}^2\\
        \overline{d_{3m+1}^\+  d_{3m+2}^\+ d_{3m+3}^\+}
    \end{array}\right)=
    \left(
    \begin{array}
      {lll}
      0&0&0\\
      1&0&1\\
      0&1&0
    \end{array}\right).
  \end{equation}
  Suppose (\ref{eq:extend-2}) holds. We will prove in this case that ${d_1^2\ldots d^2_{3m+3}}=(010)^{m+1}$.

  If $m=1$, then {$d_1^1\ldots d_6^1=101000$. By (\ref{eq:sum-one-digits}), (\ref{eq:nov-19}) and (\ref{eq:extend-2})} it follows that
  {\[
  \left(\begin{array}
    {l}
    d_1^1\ldots d_6^1\\
    d_1^2\ldots d_6^2\\
    \overline{d_1^\+\ldots d_6^\+}
  \end{array}\right)=\left(
  \begin{array}
    {ll}
    101&000\\
    010&010\\
    000&101
  \end{array}\right),
  \]
which yields $d_1^2\ldots d_6^2=(010)^2$ as desired. If $m=2$, then $d_1^1\ldots d_9^1=101001000$. By {Proposition \ref{prop:key-prop}} we can deduce that $d_1^2\ldots d_9^2=(010)^3$.}

   In the following we assume $m\ge 3$. Note that $d_{3m-2}^1d_{3m-1}^1d_{3m}^1=001$. Then by  (\ref{eq:sum-one-digits}), (\ref{eq:nov-19}) and (\ref{eq:extend-2})  we must have
 {\begin{equation}\label{eq:nov-19-1}
  \left(\begin{array}
    {l}
    d_{3m-2}^1\ldots d_{3m+3}^1\\
    d_{3m-2}^2\ldots d_{3m+3}^2\\
    \overline{d_{3m-2}^\+\ldots d_{3m+3}^\+}
  \end{array}\right)=\left(
  \begin{array}
    {ll}
    001&000\\
    010&010\\
    100&101
  \end{array}\right).
  \end{equation} So, $d_{3m-2}^2\ldots d_{3m+3}^2=(010)^2$.}
%
  Similarly, note that $d_{3m-5}^1d_{3m-4}^1d_{3m-3}^1=001$. Then we have $d_{3m-5}^2 d_{3m-4}^2 d_{3m-3}^2\in\set{010, 100}$. If $d_{3m-5}^2 d_{3m-4}^2 d_{3m-3}^2=100$, then $\overline{d_{3m-5}^\+d_{3m-4}^\+d_{3m-3}^\+}=010$, which implies that
  \[
  \overline{d_{3m-4}^\+\ldots d_{3m+3}^\+}=10100101\succ {a_1\ldots a_8},
  \]
  again leading to a contradiction with (\ref{eq:admissible-condition}). This proves $d_{3m-5}^2d_{3m-4}^2d_{3m-3}^2=010$.

  Proceeding this argument we can deduce that $d_1^2\ldots d_{3m+3}^2=(010)^{m+1}$. Symmetrically, if (\ref{eq:extend-3}) holds, then we can prove that $\overline{d_1^\+\ldots d_{3m+3}^\+}=(010)^{m+1}$.
\end{proof}

\begin{proof}
  [Proof of Theorem \ref{th:main result} (iv)]
  In view of Lemma \ref{lem:below-3}  we assume $\ell\in\N_{\ge 4}$. Suppose $\de(\beta_{3\ell+2})=(a_1\ldots a_{3\ell+2})^\f$. Then {$a_1\ldots a_{3\ell+2}$ is admissible, and  there exists a block $d_1\ldots d_{3\ell+2} \in \Omega^{*}$} satisfying (\ref{eq:admissible-condition}) and $d_1^1\ldots d_{3\ell+2}^1=a_1\ldots a_{3\ell+2}$. Note by Proposition \ref{prop:above eta-9}   that
  \begin{equation}
    \label{eq:th3n+2-1}
    (101001000)^\f\prec (a_1\ldots a_{3\ell+2})^\f\lle (101(001)^{\ell-1}00)^\f.
  \end{equation}
  {This gives $a_1\ldots a_8=10100100$.}

  Suppose on the contrary that $(a_1\ldots a_{3\ell+2})^\f\ne (101(001)^{\ell-1}00)^\f$. Then by {(\ref{eq:th3n+2-1})   there exists $m\in\set{2, 3,\ldots, \ell-1}$} such that
  \begin{equation}
    \label{eq:th3n+2-3}
    a_1\ldots a_{3m+3}=101(001)^{m-1}000.
  \end{equation}
{By} (\ref{eq:sum-one-digits}), (\ref{eq:th3n+2-3}) and Lemma \ref{lem:extended-admissible} it follows that
  \begin{equation}
    \label{eq:th3n+2-4}
     \left(
    \begin{array}
      {c}
      d_{1}^1\ldots d_{3m+3}^1\\
       d_{1}^2\ldots d_{3m+3}^2\\
        \overline{d_{1}^\+ \ldots d_{3m+3}^\+}
    \end{array}\right)=
    \left(
    \begin{array}
      {lll}
      101&(001)^{m-1}&000\\
      010&(010)^{m-1}&010\\
      000&(100)^{m-1}&101
    \end{array}\right)
  \end{equation}
  or
  \begin{equation*}
          \left(
    \begin{array}
      {c}
      d_{1}^1\ldots d_{3m+3}^1\\
       d_{1}^2\ldots d_{3m+3}^2\\
        \overline{d_{1}^\+ \ldots d_{3m+3}^\+}
    \end{array}\right)=
    \left(
    \begin{array}
      {lll}
      101&(001)^{m-1}&000\\
      000&(100)^{m-1}&101\\
       010&(010)^{m-1}&010
      \end{array}\right).
  \end{equation*}
  Without loss of generality we may assume (\ref{eq:th3n+2-4}) holds, since the other case can be proved similarly. Note by (\ref{eq:admissible-condition}) that $a_{3\ell+1}a_{3\ell+2}=00$.
  If $m=\ell-1$, then by (\ref{eq:th3n+2-4}) it follows that
  \[
  \left(
    \begin{array}
      {c}
      d_{1}^1\ldots d_{3\ell+2}^1\\
       d_{1}^2\ldots d_{3\ell+2}^2\\
        \overline{d_{1}^\+ \ldots d_{3\ell+2}^\+}
    \end{array}\right)=
    \left(
    \begin{array}
      {llll}
      101&(001)^{\ell-2}&000&00\\
      010&(010)^{\ell-2}&010&**\\
      000&(100)^{\ell-2}&101&**
    \end{array}\right)
  \]
  So, by (\ref{eq:sum-one-digits}) it follows that $\overline{d_{3\ell-2}^\+\ldots d_{3\ell+2}^\+}\lge 10101\succ a_1\ldots a_5$, leading to a contradiction with (\ref{eq:admissible-condition}).

  In the following we assume $m<\ell-1$. Note that $d_{3\ell+1}^1d_{3\ell+2}^1=a_{3\ell+1}a_{3\ell+2}=00$. Then $d_{3\ell+1}^2d_{3\ell+2}^2\in\set{01, 10}$. If it takes $01$, then by (\ref{eq:th3n+2-4}) and (\ref{eq:th3n+2-3}) we have
  \[
  d_{3\ell+2}^2d_1^2\ldots d_{3m+2}^2=101(001)^m\succ a_1\ldots a_{3m+3},
  \]
  leading to a contradiction with (\ref{eq:admissible-condition}). So, $d_{3\ell+1}^2d_{3\ell+2}^2=10$, and $\overline{d_{3\ell+1}^\+d_{3\ell+2}^\+}=01$.

  {\bf Claim.} For any ${m}\le i<\ell$ we have
  \begin{equation}
    \label{eq:th3n+2-5}
    \left(
    \begin{array}
      {c}
      d_{3i+1}^1d_{3i+2}^1 d_{3i+3}^1\\
       d_{3i+1}^2d_{3i+2}^2 d_{3i+3}^2\\
        \overline{d_{3i+1}^\+d_{3i+2}^\+ d_{3i+3}^\+}
    \end{array}\right)=
    \left(
    \begin{array}
      {c}
     010\\
     100\\
     001
    \end{array}\right)\quad\textrm{or}\quad
       \left(
    \begin{array}
      {c}
      d_{3i+1}^1d_{3i+2}^1 d_{3i+3}^1\\
       d_{3i+1}^2d_{3i+2}^2 d_{3i+3}^2\\
        \overline{d_{3i+1}^\+d_{3i+2}^\+ d_{3i+3}^\+}
    \end{array}\right)=\left(
    \begin{array}
      {c}
     001\\
     100\\
     010
    \end{array}\right).
  \end{equation}

  First we prove (\ref{eq:th3n+2-5}) for  $i=\ell-1$. Since $d_{3\ell+1}^2=1$ and $11$ is forbidden in $(d_1^2\ldots d_{3\ell+2}^2)^\f$, it follows that
  \[d_{3\ell-2}^2d_{3\ell-1}^2d_{3\ell}^2\in\set{000, 010, 100}.\]
  If $d_{3\ell-2}^2d_{3\ell-1}^2d_{3\ell}^2=010$, then by (\ref{eq:th3n+2-3}) and (\ref{eq:th3n+2-4}) we have
  \[
  d_{3\ell-1}^2\ldots d_{3\ell+2}^2d_1^2\ldots d_{3m-1}^2=101(001)^m\succ a_1\ldots a_{3m+3},
  \]
  leading to a contradiction with (\ref{eq:admissible-condition}). If $d_{3\ell-2}^2d_{3\ell-1}^2d_{3\ell}^2=000$, then either $d_{3\ell-2}^1d_{3\ell-1}^1d_{3\ell}^1=101$ which gives $d_{3\ell-2}^1\ldots d_{3\ell+2}^1d_1^1d_2^1d_3^1=10100101\succ a_1\ldots a_8$, or $\overline{d_{3\ell-2}^\+d_{3\ell-1}^\+d_{3\ell}^\+}=101$ which yields $\overline{d_{3\ell-2}^\+\ldots d_{3\ell+2}^\+}=10101\succ a_1\ldots a_5$. Both cases lead to a contradiction with (\ref{eq:admissible-condition}). Therefore, we must have $d_{3\ell-2}^2d_{3\ell-1}^2d_{3\ell}^2=100$. This together with (\ref{eq:sum-one-digits}) proves (\ref{eq:th3n+2-5}) for $i=\ell-1$.

  Next we prove (\ref{eq:th3n+2-5}) for $i=\ell-2$. Since $d_{3\ell-2}^2=1$ and $11$ is forbidden in $(d_1^2\ldots d_{3\ell+2}^2)^\f$, it follows that
  \[
  d_{3\ell-5}^2d_{3\ell-4}^2d_{3\ell-3}^2\in\set{000, 010, 100}.
  \]
  If $ d_{3\ell-5}^2d_{3\ell-4}^2d_{3\ell-3}^2=010$, then by (\ref{eq:th3n+2-3}) and (\ref{eq:th3n+2-4}) it follows that
  \[
  d_{3\ell-4}^2\ldots d_{3\ell+2}^2d_1^2\ldots d_{3m-4}^2=101(001)^m\succ a_1\ldots a_{3m+3},
  \]
  leading to a contradiction with (\ref{eq:admissible-condition}). If $d_{3\ell-5}^2d_{3\ell-4}^2d_{3\ell-3}^2=000$, then  we consider two cases in terms of (\ref{eq:th3n+2-5}) for $i=\ell-1$.
  \begin{itemize}
    \item Suppose
    \[
    \left(
    \begin{array}
      {c}
      d_{3\ell-2}^1\ldots d_{3\ell+2}^1\\
       d_{3\ell-2}^2\ldots d_{3\ell+2}^2\\
        \overline{d_{3\ell-2}^\+\ldots d_{3\ell+2}^\+}
    \end{array}\right)=
    \left(
    \begin{array}
      {c}
     01000\\
     10010\\
     00101
    \end{array}\right).
    \]
    Then by using $d_{3\ell-5}^2d_{3\ell-4}^2d_{3\ell-3}^2=000$ it follows that either $d_{3\ell-5}^1d_{3\ell-4}^1d_{3\ell-3}^1=101$ which gives $d_{3\ell-5}^1\ldots d_{3\ell-1}^1=10101\succ a_1\ldots a_5$, or $\overline{d_{3\ell-5}^\+d_{3\ell-4}^\+d_{3\ell-3}^\+}=101$ which implies
    $\overline{d_{3\ell-5}^\+\ldots d_{3\ell+2}^\+}=10100101\succ a_1\ldots a_8$. Both cases  lead to a contradiction with (\ref{eq:admissible-condition}).

    \item Suppose
      \[
    \left(
    \begin{array}
      {c}
      d_{3\ell-2}^1\ldots d_{3\ell+2}^1\\
       d_{3\ell-2}^2\ldots d_{3\ell+2}^2\\
        \overline{d_{3\ell-2}^\+\ldots d_{3\ell+2}^\+}
    \end{array}\right)=
    \left(
    \begin{array}
      {c}
     00100\\
     10010\\
     01001
    \end{array}\right).
    \]
    Then by using $d_{3\ell-5}^2d_{3\ell-4}^2d_{3\ell-3}^2=000$ it follows that either $d_{3\ell-5}^1d_{3\ell-4}^1d_{3\ell-3}^1=101$ which implies
    $d_{3\ell-5}^1\ldots d_{3\ell+2}^1 d_1^1d_2^1d_3^1 =10100100101\succ a_1\ldots a_{11},$
    or $\overline{d_{3\ell-5}^\+d_{3\ell-4}^\+d_{3\ell-3}^\+}=101$ which gives $\overline{d_{3\ell-5}^\+\ldots d_{3\ell-1}^\+}=10101\succ a_1\ldots a_5$. Both cases again lead to a contradiction with (\ref{eq:admissible-condition}).
  \end{itemize}
  Therefore, $d_{3\ell-5}^2d_{3\ell-4}^2d_{3\ell-3}^2=100$, and from this we deduce that (\ref{eq:th3n+2-5}) holds for $i=\ell-2$.

  Proceeding this argument for $i=\ell-3, \ell-4,\ldots, {m+1}$ we can prove (\ref{eq:th3n+2-5})   for all ${m+1}\le i<\ell$. This establishes the claim. But then by the claim and (\ref{eq:th3n+2-4}) {it follows that $d_1^2\ldots d_{3\ell+2}^2=(010)^{m+1}(100)^{\ell-m-1}10$. This together with (\ref{eq:th3n+2-3}) implies that}
  \[
  d_{3m+2}^2\ldots d_{3\ell+2}^2d_1^2\ldots d_{3m-1}^2=101(001)^{\ell-1}\succ a_1\ldots a_{3\ell},
  \]
  again leading to a contradiction with (\ref{eq:admissible-condition}). Hence, $(a_1\ldots a_{3\ell+2})^\f=(101(001)^{\ell-1}00)^\f$, completing the proof.
\end{proof}

\section{Minimal base for $\U_\beta$ to have a sequence of period  $3\N+1$}\label{sec:3N+1}
In this section we will consider $\beta_{3\ell+1}$ for any $\ell\in\N$, and prove Theorem \ref{th:main result} (iii).
First we prove it for $\ell\in\set{1,2}$.
\begin{lemma}
  \label{lem:3n+1bounded-4}
  If $\ell\in\set{1,2}$, then $\de(\beta_{3\ell+1})=(101(001)^{\lfloor\frac{\ell-1}{2}\rfloor}(010)^{\lceil\frac{\ell-1}{2}\rceil} 0)^\f.$
\end{lemma}
\begin{proof}
If $\ell=1$, then we let $\de(\beta_4)=(a_1\ldots a_4)^\f$. By   Theorem \ref{th:lower-higher-orders} it follows that
\[
1010\ldots=\de(\beta_c)\lle (a_1\ldots a_4)^\f\lle \de(\beta_2)=(10)^\f,
\]
{which implies} $\de(\beta_4)=(10)^\f$ as desired.

 Next we consider $\ell=2$.    Let $\de(\beta_7)=(a_1\ldots a_7)^\f$. By  Theorem \ref{th:lower-higher-orders} we obtain that
    \[
    101001000\ldots=\de(\beta_c)\lle (a_1\ldots a_7)^\f\lle \de(\beta_2)=(10)^\f.
    \]
     This implies $a_1 a_2a_3=101$ and $a_4a_5a_6\in\set{001, 010}$. If we take $a_4a_5a_6=001$, then by   (\ref{eq:admissible-condition}) we must have $a_7=0$, and thus $a_6a_7a_1a_2a_3=10101\succ a_1\ldots a_5$, again leading to a contradiction with (\ref{eq:admissible-condition}). So, $a_1\ldots a_7=1010100$ as desired.
\end{proof}
In the following we assume $\ell\ge 3$.
In the following proposition we give the bounds on $\beta_{3\ell+1}$   for     $\ell\in\N$.
\begin{proposition}
  \label{prop:above eta-9+}\mbox{}
   For any $\ell\in\N$, we have
    \begin{equation}\label{eq:3l+1}
    (101001000)^\f=\de(\beta_9)\prec\de(\beta_{3\ell+1})\lle (101(001)^{\lfloor\frac{\ell-1}{2}\rfloor}(010)^{\lceil\frac{\ell-1}{2}\rceil} 0)^\f.
    \end{equation}
\end{proposition}
\begin{proof}
    By Lemma \ref{lem:3n+1bounded-4} it suffices to prove (\ref{eq:3l+1}) for $\ell\ge 3$. Take
  \[(d_1\ldots d_{3\ell+1})=(\al_1\al_2\al_1(\al_0\al_2\al_1)^{\lfloor\frac{\ell-1}{2}\rfloor}(\al_0\al_1\al_2)^{\lceil\frac{\ell-1}{2}\rceil}\al_0)^\f.\]
  Then
  \begin{align*}
  (d_1^1\ldots d_{3\ell+1}^1)^\f&=(101(001)^{\lfloor\frac{\ell-1}{2}\rfloor}(010)^{\lceil\frac{\ell-1}{2}\rceil} 0)^\f,\\
  (d_1^2\ldots d_{3\ell+1}^2)^\f&=(010(010)^{\lfloor\frac{\ell-1}{2}\rfloor}(001)^{\lceil\frac{\ell-1}{2}\rceil} 0)^\f,\\
  (\overline{d_1^\+ \ldots  d_{3\ell+1}^\+})^\f&=(000(100)^{\lfloor\frac{\ell-1}{2}\rfloor}(100)^{\lceil\frac{\ell-1}{2}\rceil}1)^\f.
  \end{align*}
  By Lemma \ref{lem:lexicographical-periodic-U} one can verify that $(d_1\ldots d_{3\ell+1})^\f\in \U_\beta$ if $\de(\beta)\succ (101(001)^{\lfloor\frac{\ell-1}{2}\rfloor}(010)^{\lceil\frac{\ell-1}{2}\rceil} 0)^\f$. This proves the second inequality of (\ref{eq:3l+1}).

 To prove the first inequality of (\ref{eq:3l+1}) we write $\delta(\beta_{3\ell+1})= (a_1\ldots a_{3\ell+1})^\f$.
  Suppose on the contrary that $\de(\beta_{3\ell+1})\lle (101001000)^\f$. Then by Theorem \ref{th:lower-higher-orders} (ii) it follows that
\[
101001000101\ldots=\de(\beta_c)\lle (a_1\ldots a_{3\ell+1})^\f\lle (101001000)^\f,
\] which implies $a_1\ldots a_9=101001000$. Note that $a_1\ldots a_{3\ell+1}$ is an admissible block, and then it has a representation $d_1\ldots d_{3\ell+1}\in \Omega^{3\ell+1}$. By {Proposition \ref{prop:key-prop}} it follows that
$
a_1\ldots a_{3\ell}\in\mathcal B^*(X),$
  and
  \begin{equation}\label{eq:3n+1-1}
  d_1^1\ldots d_{3\ell}^1=a_1\ldots a_{3\ell}, \quad d_1^2\ldots d_{3\ell}^2=(010)^\ell\quad\textrm{and}\quad \overline{d_1^\+}\ldots \overline{d_{3\ell}^\+}=\Theta(a_1\ldots a_{3\ell}).
\end{equation}
Note by (\ref{eq:admissible-condition}) that $a_{3\ell+1}=0$.
We will finish the proof by considering the following four cases: $a_{3\ell-2}a_{3\ell-1}a_{3\ell}\in\set{101, 001, 000, 100}$.

Case I. $a_{3\ell-2}a_{3\ell-1} a_{3\ell}=101$. Then $a_{3\ell-2}\ldots a_{3\ell+1}a_1=10101\succ a_1\ldots a_{5}$, leading to a contradiction with (\ref{eq:admissible-condition}).

Case II. $a_{3\ell-2}a_{3\ell-1}a_{3\ell}=001$. Then $a_{3\ell}a_{3\ell+1}a_1a_2a_3=10101\succ a_1\ldots a_5$, again leading to a contradiction with (\ref{eq:admissible-condition}).

Case III. $a_{3\ell-2}a_{3\ell-1}a_{3\ell}=000$. Then by (\ref{eq:3n+1-1}) we must have
\[
\left(\begin{array}
  {c}
  d_{3\ell-2}^1\ldots d_{3\ell+1}^1\\
  d_{3\ell-2}^2\ldots d_{3\ell+1}^2\\
  \overline{d_{3\ell-2}^\+}\ldots \overline{d_{3\ell+1}^\+}
\end{array}\right)=
\left(
\begin{array}
  {cccc}
  0&0&0&0\\
  0&1&0&1\\
  1&0&1&0
\end{array}\right),
\]
which  yields that
\[
d_{3\ell-1}^2d_{3\ell}^2d_{3\ell+1}^2d_1^2d_2^2=10101  \succ a_1\ldots a_5,
\]
  leading to a contradiction.

Case IV. $a_{3\ell-2}a_{3\ell-1}a_{3\ell}=100$. Then {$a_{3\ell-5}a_{3\ell-4} a_{3\ell-3}=000$}, and thus by (\ref{eq:3n+1-1}) we   have
\[
\left(\begin{array}
  {c}
  d_{3\ell-5}^1\ldots d_{3\ell+1}^1\\
  d_{3\ell-5}^2\ldots d_{3\ell+1}^2\\
  \overline{d_{3\ell-5}^\+}\ldots \overline{d_{3\ell+1}^\+}
\end{array}\right)=
\left(
\begin{array}
  {ccc}
  000&100&0\\
  010&010&1\\
  101&001&0
\end{array}\right).
\]
Again, this together with (\ref{eq:3n+1-1}) implies  that $d_{3\ell-1}^2d_{3\ell}^2d_{3\ell+1}^2d_1^2d_2^2=10101  \succ a_1\ldots a_5$,   leading to a contradiction.

By Cases I-IV and   Theorem \ref{th:lower-higher-orders} we conclude that $\de(\beta_{3\ell+1})\succ(101001000)^\f$ as desired.
\end{proof}
{Recall from Theorem \ref{th:eta-period-3n} that $\beta_{3\ell}\le \beta_{9}$ for all $\ell\in\N$. By Proposition \ref{prop:above eta-9+} and Lemma \ref{lem:delta-beta} it follows that $\beta_{3\ell+1}>\beta_9$ for all $\ell\in\N$.}

In the following we assume $\ell\ge 3$. Before proving Theorem \ref{th:main result} (iii) we need the following two lemmas.
\begin{lemma}
  \label{lem:after010}
  Let $a_1\ldots a_{3\ell+1}$ be an admissible {block} with $\ell\ge 3$ and
  \begin{equation}\label{eq:(3n+1)-1}
  a_1\ldots a_{3(\lfloor\frac{\ell-1}{2}\rfloor+1)}=101(001)^{\lfloor\frac{\ell-1}{2}\rfloor}.
  \end{equation} If
  $a_{3(m-j)}a_{3(m-j)+1}\ldots a_{3m}=1(010)^{j}$ for some $j\ge 1$ and $\lfloor\frac{\ell-1}{2}\rfloor+1+j\le m <\ell$, then
  \[a_{3m+1}a_{3m+2}a_{3m+3}\in\set{001, 010}.\]
\end{lemma}
{By Lemma \ref{lem:after010} it follows that if $a_1\ldots a_{3\ell+1}$ is admissible and it has a prefix $101(001)^{\lfloor\frac{\ell-1}{2}\rfloor}(010)^j$, then it can only be followed by a block in $\set{001, 010}$.}
\begin{proof}
  Note that $a_{3(m-j)}a_{3(m-j)+1}\ldots a_{3m}=1(010)^{j}$ and $j\le m-1-\lfloor\frac{\ell-1}{2}\rfloor\le \lfloor\frac{\ell-1}{2}\rfloor$. By (\ref{eq:(3n+1)-1}) and (\ref{eq:admissible-condition}) it follows that
  \[
  a_{3m+1}a_{3m+2}a_{3m+3}\in\set{000,001, 010}.
  \]
  So it suffices to prove $a_{3m+1}a_{3m+2}a_{3m+3}\ne 000$.

  Suppose on the contrary that $a_{3m+1}a_{3m+2}a_{3m+3}=000$. Since $a_1\ldots a_{3\ell+1}$ is admissible, there exists an aperiodic {block} $d_1\ldots d_{3\ell+1}\in \Omega^{3\ell+1}$ satisfying (\ref{eq:admissible-condition}) and $d_1^1\ldots d_{3\ell+1}^1=a_1\ldots a_{3\ell+1}$. Note that $a_{3m}\ldots a_{3m+3}=0000$. By (\ref{eq:admissible-condition}) and (\ref{eq:sum-one-digits}) it follows that either $d_{3m}^2\ldots d_{3m+3}^2=1010$ or $\overline{d_{3m}^\+\ldots d_{3m+3}^\+}=1010$. Since the proofs for the two cases are similar, we may assume $d_{3m}^2\ldots d_{3m+3}^2=1010$, and in this case we claim that
  \begin{equation}\label{eq:(3n+1)-2}
  d_{3m+1}^2\ldots d_{3\ell}^2=(010)^{\ell-m}.
  \end{equation}

  Observe {by (\ref{eq:sum-one-digits})} that
  \[
       \left(
    \begin{array}
      {c}
      d_{3m}^1\ldots d_{3m+6}^1\\
       d_{3m}^2\ldots d_{3m+6}^2\\
        \overline{d_{3m}^\+\ldots d_{3m+6}^\+}
    \end{array}\right)=
    \left(
    \begin{array}
      {lll}
     0&000&***\\
    1&010&*** \\
     0&101&***
    \end{array}\right).
  \]
  Note by (\ref{eq:(3n+1)-1}) and (\ref{eq:admissible-condition}) that $d_{3m+4}^2=\overline{d_{3m+4}^\+}=0$, and then $d_{3m+4}^1=1$ by (\ref{eq:sum-one-digits}). Again by (\ref{eq:(3n+1)-1}) and (\ref{eq:admissible-condition}) we have $d_{3m+5}^1=0=\overline{d_{3m+5}^\+}$. So, $d_{3m+5}^2=1$. For $d_{3m+6}$ we  have a choice   $d_{3m+6} \in\set{\al_0,\al_1}$.
  \begin{enumerate}[{\rm(i)}]
    \item If $d_{3m+6}=\al_0$, then by (\ref{eq:sum-one-digits}), (\ref{eq:admissible-condition}) and (\ref{eq:(3n+1)-1}) it follows that
    \[
       \left(
    \begin{array}
      {c}
      d_{3m}^1\ldots d_{3m+9}^1\\
       d_{3m}^2\ldots d_{3m+9}^2\\
        \overline{d_{3m}^\+\ldots d_{3m+9}^\+}
    \end{array}\right)=
    \left(
    \begin{array}
      {llll}
     0&000&100&10*\\
    1&010&010&010 \\
     0&101&001&00*
    \end{array}\right).
  \]
  This implies that $d_{3m+9}\in\set{\al_0, \al_1}$.

  \item If $d_{3m+6}=\al_1$, then by (\ref{eq:sum-one-digits}), (\ref{eq:admissible-condition}) and (\ref{eq:(3n+1)-1}) we obtain that
      \[
       \left(
    \begin{array}
      {c}
      d_{3m}^1\ldots d_{3m+9}^1\\
       d_{3m}^2\ldots d_{3m+9}^2\\
        \overline{d_{3m}^\+\ldots d_{3m+9}^\+}
    \end{array}\right)=
    \left(
    \begin{array}
      {llll}
     0&000&101&00*\\
    1&010&010&010 \\
     0&101&000&10*
    \end{array}\right).
  \]
  Again this gives that $d_{3m+9}\in\set{\al_0, \al_1}$.
  \end{enumerate}

  Proceeding this argument we conclude that $d_{3m+1}^2\ldots d_{3\ell}^2=(010)^{\ell-m}$, establishing (\ref{eq:(3n+1)-2}). So,
  \[
  d_{3m}^2\ldots d_{3\ell}^2=101(001)^{\ell-m-1}0.
  \]
  Note that $0\le \ell-m-1<\lfloor\frac{\ell-1}{2}\rfloor$. By (\ref{eq:(3n+1)-1}) and (\ref{eq:admissible-condition}) it follows that $d_{3\ell}^2d_{3\ell+1}^2=00$. However, by (\ref{eq:(3n+1)-1}) and (\ref{eq:admissible-condition}) we must have $d_{3\ell}^1d_{3\ell+1}^1{=a_{3\ell}a_{3\ell+1}}=00$. Thus, by (\ref{eq:sum-one-digits}) it gives that
  $\overline{d_{3\ell}^\+d_{3\ell+1}^\+}=11$, leading to a contradiction with (\ref{eq:admissible-condition}).

  Therefore, $a_{3m+1}a_{3m+2}a_{3m+3}\ne 000$, completing the proof.
  \end{proof}

\begin{lemma}
  \label{lem:after001}
  Let $a_1\ldots a_{3\ell+1}$ be an admissible {block} with $\ell\ge 3$ and
  \begin{equation}\label{eq:(3n+1)-3}
  a_1\ldots a_{3(\lfloor\frac{\ell-1}{2}\rfloor+1)}=101(001)^{\lfloor\frac{\ell-1}{2}\rfloor}.
  \end{equation} If
 $a_{3(m-j)}a_{3(m-j)+1}\ldots a_{3m}=0(001)^{j}$ for some $j\ge 1$ and $\lfloor\frac{\ell-1}{2}\rfloor+{2}+j\le m <\ell$, then
  \[a_{3m+1}a_{3m+2}a_{3m+3}\in\set{001, 010}.\]
\end{lemma}
{By Lemmas \ref{lem:after010} and \ref{lem:after001} it follows that if $a_1\ldots a_{3\ell+1}$ is admissible, and it has a prefix $101(001)^{\lfloor\frac{\ell-1}{2}\rfloor}(010)^{j_1}(001)^{j_2}$ with $j_1, j_2\ge 1$, then it can only be followed by a block in $\set{001, 010}$.}

\begin{proof}
  Since $a_{3m}=1$, by (\ref{eq:(3n+1)-3}) and (\ref{eq:admissible-condition}) it follows that $a_{3m+1}a_{3m+2}a_{3m+3}\in\set{000, 001, 010}$. In the following it suffices to prove $a_{3m+1}a_{3m+2}a_{3m+3}\ne 000$.

  Suppose on the contrary that $a_{3m+1}a_{3m+2}a_{3m+3}=000$. Since $a_1\ldots a_{3\ell+1}$ is admissible, there exists an aperiodic {block} $d_1\ldots d_{3\ell+1}$ satisfying (\ref{eq:admissible-condition}) and $d_1^1\ldots d_{3\ell+1}^1=a_1\ldots a_{3\ell+1}$. Note that $a_{3(m-j)}\ldots a_{3m}=0(001)^j$. By (\ref{eq:sum-one-digits}) and (\ref{eq:admissible-condition}) it follows that either
  \[d_{3(m-j)}^2d_{3(m-j)+1}^2d_{3(m-j)+2}^2=101\quad\textrm{or}\quad \overline{d_{3(m-j)}^\+d_{3(m-j)+1}^\+d_{3(m-j)+2}^\+}=101.\]
 By symmetry we may assume $d_{3(m-j)}^2d_{3(m-j)+1}^2d_{3(m-j)+2}^2=101$.\\
 \medskip

 {\bf Claim.}
 \[
 d_{3(m-j)+1}^2\ldots d_{3\ell}^2=(010)^{\ell-m+j}.
 \]

 Observe {by (\ref{eq:sum-one-digits})} that
       \[
       \left(
    \begin{array}
      {c}
      d_{3(m-j)}^1\ldots d_{3m}^1\\
       d_{3(m-j)}^2\ldots d_{3m}^2\\
        \overline{d_{3(m-j)}^\+\ldots d_{3m}^\+}
    \end{array}\right)=
    \left(
    \begin{array}
      {lll}
     0&001&(001)^{j-1}\\
    1&010&*** \\
     0&100&***
    \end{array}\right).
  \]
 Then by (\ref{eq:sum-one-digits}), (\ref{eq:admissible-condition}) and (\ref{eq:(3n+1)-3}) it follows that
  \[
  d_{3(m-j)+4}d_{3(m-j)+5}d_{3(m-j)+6}=\al_0\al_2\al_1.
  \]
{Note that $j<m-\lfloor\frac{\ell-1}{2}\rfloor\le \lfloor\frac{\ell-1}{2}\rfloor+1$.} Proceeding using (\ref{eq:sum-one-digits}), (\ref{eq:admissible-condition}) and (\ref{eq:(3n+1)-3}) we can deduce that
 $d_{3(m-j)+1}\ldots d_{3m}=(\al_0\al_2\al_1)^j$.

 Note that $d_{3m+1}^1d_{3m+2}^1d_{3m+3}^1=a_{3m+1}a_{3m+2}a_{3m+3}=000$. By (\ref{eq:sum-one-digits}) and (\ref{eq:admissible-condition}) it follows that
         \[
       \left(
    \begin{array}
      {c}
      d_{3(m-j)}^1\ldots d_{3m+3}^1\\
       d_{3(m-j)}^2\ldots d_{3m+3}^2\\
        \overline{d_{3(m-j)}^\+\ldots d_{3m+3}^\+}
    \end{array}\right)=
    \left(
    \begin{array}
      {lll}
     0&(001)^j&000\\
    1&(010)^j&010 \\
     0&(100)^j&101
    \end{array}\right).
  \]
  Again by (\ref{eq:sum-one-digits}), (\ref{eq:admissible-condition}) and (\ref{eq:(3n+1)-3})  we obtain that $d_{3m+4}d_{3m+5}=\al_1\al_2$ and $d_{3m+6}\in\set{\al_0, \al_1}$.
  \begin{enumerate}
    [{\rm(i)}]
    \item If $d_{3m+6}=\al_0$, then by (\ref{eq:sum-one-digits}), (\ref{eq:admissible-condition}) and (\ref{eq:(3n+1)-3})  it follows that
             \[
       \left(
    \begin{array}
      {c}
      d_{3(m-j)}^1\ldots d_{3m+9}^1\\
       d_{3(m-j)}^2\ldots d_{3m+9}^2\\
        \overline{d_{3(m-j)}^\+\ldots d_{3m+9}^\+}
    \end{array}\right)=
    \left(
    \begin{array}
      {lllll}
     0&(001)^j&000&100&10*\\
    1&(010)^j&010 &010&010\\
     0&(100)^j&101&001&00*
    \end{array}\right),
  \]
which implies  $d_{3m+9}\in\set{\al_0, \al_1}$.

  \item If $d_{3m+6}=\al_1$, then by (\ref{eq:sum-one-digits}), (\ref{eq:admissible-condition}) and (\ref{eq:(3n+1)-3}) we can deduce that
            \[
       \left(
    \begin{array}
      {c}
      d_{3(m-j)}^1\ldots d_{3m+9}^1\\
       d_{3(m-j)}^2\ldots d_{3m+9}^2\\
        \overline{d_{3(m-j)}^\+\ldots d_{3m+9}^\+}
    \end{array}\right)=
    \left(
    \begin{array}
      {lllll}
     0&(001)^j&000&101&00*\\
    1&(010)^j&010 &010&010\\
     0&(100)^j&101&000&10*
    \end{array}\right),
  \]
which yields that $d_{3m+9}\in\set{\al_0, \al_1}$.
  \end{enumerate}

  Iterating the above arguments we conclude that $d_{3(m-j)+1}^2\ldots d_{3\ell}^2=(010)^{\ell-m+j}$, establishing the claim. Therefore,
  \[
  d_{3(m-j)}^2\ldots d_{3\ell}^2=101(001)^{\ell-m+j-1}0.
  \]
 {Since $m-j\ge \lfloor\frac{\ell-1}{2}\rfloor+2$, we have} $0\le \ell-m+j-1<\lfloor\frac{\ell-1}{2}\rfloor$. Then by (\ref{eq:sum-one-digits}), (\ref{eq:admissible-condition}) and (\ref{eq:(3n+1)-3}) it follows that $d_{3\ell}^2d_{3\ell+1}^2=00$. However, since $d_{3\ell}^1d_{3\ell+1}^1{=a_{3\ell}a_{3\ell+1}}=00$, we must have $\overline{d_{3\ell}^\+d_{3\ell+1}^\+}=11$, leading  to a contradiction with (\ref{eq:admissible-condition}).
    Therefore, $a_{3m+1}a_{3m+2}a_{3m+3}\ne 000$ as desired.
\end{proof}

\begin{proof}
  [Proof of Theorem \ref{th:main result} (iii)]
By Lemma \ref{lem:3n+1bounded-4} it suffices to consider $\ell\ge 3$. Let $\de(\beta_{3\ell+1})=(a_1\ldots a_{3\ell+1})^\f$. Then $a_1\ldots a_{3\ell+1}$ is admissible, and  there exists an aperiodic {block} $d_1\ldots d_{3\ell+1}\in \Omega^{3\ell+1}$ satisfying (\ref{eq:admissible-condition}) and $d_1^1\ldots d_{3\ell+1}^1=a_1\ldots a_{3\ell+1}$. By Proposition \ref{prop:above eta-9+} it follows that
\begin{equation}
  \label{eq:th3n+1-1}
  (101001000)^\f\prec (a_1\ldots a_{3\ell+1})^\f\lle (101(001)^{\lfloor\frac{\ell-1}{2}\rfloor}(010)^{\lceil\frac{\ell-1}{2}\rceil}0)^\f.
\end{equation}
 Then by (\ref{eq:th3n+1-1}) there exists $m\in\set{2,3,\ldots, \ell-1}$ such that
$a_1\ldots a_{3m+3}=101(001)^{m-1}000,$
or there exists $m\in\set{\lfloor\frac{\ell-1}{2}\rfloor+1, \ldots, \ell-1}$ such that
$
a_1\ldots a_{3m+3}=101(001)^{m-1}010.
$
In terms of the two possibilities we split the proof into two cases.

Case I. There exists $m\in\set{2,3,\ldots, \ell-1}$ such that
$a_1\ldots a_{3m+3}=101(001)^{m-1}000$. By Lemma \ref{lem:extended-admissible} it follows that
\begin{equation}\label{eq:th3n+1-2}
d_1^2\ldots d_{3m+3}^2=(010)^{m+1}\quad\textrm{and}\quad \overline{d_1^\+\ldots d_{3m+3}^\+}=000(100)^{m-1}101,
\end{equation} or
\[
d_1^2\ldots d_{3m+3}^2=000(100)^{m-1}101\quad\textrm{and}\quad \overline{d_1^\+\ldots d_{3m+3}^\+}=(010)^{m+1}.
\]
Without loss of generality we may assume
(\ref{eq:th3n+1-2}) holds.
Since $a_1\ldots a_{3\ell+1}$ is admissible, by (\ref{eq:admissible-condition}) we must have $a_{3\ell}a_{3\ell+1}=00$, and then $d_{3\ell}^2d_{3\ell+1}^2\in\set{01, 10}$.
If $d_{3\ell}^2d_{3\ell+1}^2=01$, then by (\ref{eq:th3n+1-2}) we have
\[
d_{3\ell+1}^2d_1^2\ldots d_{3m+2}^2=101(001)^m\succ a_1\ldots a_{3m+3},
\]
leading to a contradiction with (\ref{eq:admissible-condition}). So, $d_{3\ell}^2d_{3\ell+1}^2=10$. Then by the same argument as in the proof of (\ref{eq:th3n+2-5}) one can deduce that
\[
d_{3m+4}^2\ldots d_{3\ell+1}^2=(001)^{\ell-m-1}0,
\]
which implies that $\overline{d_{3m+4}^\+d_{3m+5}^\+}\in\set{01, 10}$. Then by (\ref{eq:th3n+1-2}) it follows that
\[\overline{d_{3m+1}^\+\ldots d_{3m+5}^\+}\lge 10101\succ a_1\ldots a_5,\] leading to a contradiction with (\ref{eq:admissible-condition}).

Case II. There exists $m\in\set{\lfloor\frac{\ell-1}{2}\rfloor+1, \ldots, \ell-1}$ such that
$
a_1\ldots a_{3m+3}=101(001)^{m-1}010.
$
Then $a_1\ldots a_{3(\lfloor\frac{\ell-1}{2}\rfloor+1)}=101(001)^{\lfloor\frac{\ell-1}{2}\rfloor}$ and $a_{3m}\ldots a_{3m+3}=1010$ with $m\ge \lfloor\frac{\ell-1}{2}\rfloor+1$. By Lemma \ref{lem:after010} it follows that
\[
a_{3m+4}a_{3m+5}a_{3m+6}\in\set{001, 010}.
\]
If $a_{3m+4}a_{3m+5}a_{3m+6}=010$, then we can apply Lemma \ref{lem:after010} again to get $a_{3m+7}a_{3m+8}a_{3m+9}\in\set{001, 010}$. If $a_{3m+4}a_{3m+5}a_{3m+6}=001$, then $a_{3m+3}\ldots a_{3m+6}=0001$. By Lemma \ref{lem:after001} we can deduce that $a_{3m+7}a_{3m+8}a_{3m+9}\in\set{001, 010}$.

Proceeding this argument with applications of Lemmas \ref{lem:after010} and \ref{lem:after001} we obtain that
\[
a_{3m+4}\ldots a_{3\ell}\in\set{001, 010}^{\ell-m-1}.
\]
Suppose $j$ is the largest index in which $a_{3j+1}a_{3j+2}a_{3j+3}=001$. Then by using $a_{3\ell}a_{3\ell+1}=00$ it follows that $m-1\le j\le \ell-2$ and
\begin{equation}\label{eq:th3n+1-3}
a_{3j+1}\ldots a_{3\ell+1}=001(010)^{\ell-j-1}0.
\end{equation}
This implies that
\begin{equation}\label{eq:th3n+1-4}
a_{3j+3}\ldots a_{3\ell+1}a_1a_2a_3=101(001)^{\ell-j-1}01.
\end{equation}
Note that $m-1\le j\le\ell-2$ and $\lfloor\frac{\ell-1}{2}\rfloor+1\le m\le \ell-1$. Then
\begin{equation}\label{eq:oct-17}
1\le \ell-j-1\le \ell-m\le \ell-1-\lfloor\frac{\ell-1}{2}\rfloor=\lceil\frac{\ell-1}{2}\rceil,
\end{equation}
which gives $m\ge \lfloor\frac{\ell-1}{2}\rfloor+1$. In the following we consider two subcases.

Case (IIa). $m\ge \lfloor\frac{\ell-1}{2}\rfloor+2$. Then  by (\ref{eq:oct-17}) it follows that
\[
\ell-j-1\le \ell-m\le \ell-\lfloor\frac{\ell-1}{2}\rfloor-2=\lceil\frac{\ell-1}{2}\rceil-1<\lfloor\frac{\ell-1}{2}\rfloor+1\le m-1.
\]
By (\ref{eq:th3n+1-4}) and using $a_1\ldots a_{3m+3}=101(001)^{m-1}010$ we obtain that
\[
a_{3j+3}\ldots a_{3\ell+1}a_1a_2a_3 \succ a_1\ldots a_{3(\ell-j)+2},
\]
  leading to a contradiction with (\ref{eq:admissible-condition}).

  Case (IIb). $m=\lfloor\frac{\ell-1}{2}\rfloor+1$. Then by (\ref{eq:oct-17}) we have
  \[
  \ell-j-1\le \ell-m=\ell-\lfloor\frac{\ell-1}{2}\rfloor-1=\lceil\frac{\ell-1}{2}\rceil.
  \]
  If $\ell-j-1=\lceil\frac{\ell-1}{2}\rceil$, then by  (\ref{eq:th3n+1-3}) and $a_1\ldots a_{3m+3}=101(001)^{m-1}010$ it follows that
  \[a_1\ldots a_{3\ell+1}=101(001)^{\lfloor\frac{\ell-1}{2}\rfloor}(010)^{\lceil\frac{\ell-1}{2}\rceil}0 \] as desired.
Otherwise,
\begin{equation}
  \label{eq:nov-23}
  \ell-j-1\le\lceil\frac{\ell-1}{2}\rceil-1\le\lfloor\frac{\ell-1}{2}\rfloor=m-1.
\end{equation}
If $\ell-j-1<m-1$, then by the same argument as in Case (IIa) one can deduce a contradiction. If $\ell-j-1=m-1$, then by (\ref{eq:nov-23}) we have $\ell=2m$ and $j=m$. So by (\ref{eq:th3n+1-3}) and using $a_1\ldots a_{3m+3}=101(001)^{m-1}010$ it follows that
\[
010=a_{3m+1}a_{3m+2}a_{3m+3}=001,
\]
again leading to a contradiction.

Therefore, by Cases I and II we obtain that $a_1\ldots a_{3\ell+1}=101(001)^{\lfloor\frac{\ell-1}{2}\rfloor}(010)^{\lceil\frac{\ell-1}{2}\rceil}0$, completing the proof.
\end{proof}

\begin{proof}
  [Proof of Proposition \ref{prop:asymptotic-beta-k}]
 First we prove (i).  Note by Theorem \ref{th:main result} (ii) that for $m\in\N$ and $n\in\N_0$, we have
$\de(\beta_{3(2m+1)2^n})=(\t_{n+2}^+\Theta(\t_{n+1}^+)\t_{n+2}^{m-1})^\f$. Then   for any $m\in\N$,
\[
\de(\beta_{3(2m+3)2^n})=(\t_{n+2}^+\Theta(\t_{n+1}^+)\t_{n+2}^m)^\f\prec(\t_{n+2}^+\Theta(\t_{n+1}^+)\t_{n+2}^{m-1})^\f=\de(\beta_{3(2m+1)2^n}),
\]
and by (\ref{eq:def-tn}) we obtain that
\[
\de(\beta_{3(2m+1)2^n})\to\t_{n+2}^+\Theta(\t_{n+1}^+)\t_{n+2}^\f=\t_{n+2}^+\Theta(\t_{n+1}^+)(\t_{n+1}^+\Theta(\t_{n+1}^+))^\f=\t_{n+2}^+\Theta(\t_{n+2})^\f=\de(\hat\beta_n)
\] {as $m\to\f$.}
So, by Lemma \ref{lem:delta-beta} it follows that $\beta_{3(2m+1)2^n}\searrow \hat\beta_n$ as $\N\ni m\to\f$.
Similarly, by (\ref{eq:def-tn}) we have
\[
\de(\hat\beta_{n+1})=\t_{n+3}^+\Theta(\t_{n+3})^\f=\t_{n+2}^+\Theta(\t_{n+2}^+)\Theta(\t_{n+3})^\f\prec \t_{n+2}^+\Theta(\t_{n+2})^\f=\de(\hat\beta_n),
\]
and by (\ref{eq:lambda-i}) we have $\de(\hat\beta_n)\to\la_1\la_2\ldots =\de(\beta_c)$ as $n\to\f$. Again by Lemma \ref{lem:delta-beta} we deduce that $\hat\beta_n\searrow \beta_c$ as $n\to\f$.

 Next we consider (ii). If $m=0$, then by Theorem \ref{th:main result} (ii) {and (\ref{eq:def-tn})} it follows that
 \[
 \de(\beta_{3\cdot 2^n})=\t_{n+1}^\f\prec (\t_{n+1}^+\Theta(\t_{n+1}^+))^\f=\t_{n+2}^\f=\de(\beta_{3\cdot 2^{n+1}}),
 \]
 and by (\ref{eq:lambda-i}) we have $\de(\beta_{3\cdot 2^n})=\t_{n+1}^\f\to\de(\beta_c)$ as $n\to\f$. This together with Lemma \ref{lem:delta-beta} implies that $\beta_{3\cdot 2^n}\nearrow\beta_c$ as $n\to\f$. If $m\in\N$, then by Theorem \ref{th:main result} (ii) and (\ref{eq:def-tn}) it follows that
 \begin{align*}
   \de(\beta_{3(2m+1)2^{n+1}})&=(\t_{n+3}^+\Theta(\t_{n+2}^+)\t_{n+3}^{m-1})^\f =(\t_{n+2}^+\Theta(\t_{n+1}^+)\t_{n+1}^+\,\Theta(\t_{n+2}^+)\t_{n+3}^{m-1})^\f\\
   &\succ(\t_{n+2}^+\Theta(\t_{n+1}^+)\t_{n+2}^{m-1})^\f=\de(\beta_{3(2m+1)2^n}),
 \end{align*}
   and $\de(\beta_{3(3m+1)2^n})\to\de(\beta_c)$ as $n\to\f$ by (\ref{eq:lambda-i}). So, by Lemma \ref{lem:delta-beta} we conclude that $\beta_{3(2m+1)2^n}\searrow\beta_c$ as $n\to\f$.

   Finally we prove (iii). Note by Theorem \ref{th:main result} (iii) and (iv) that
   \begin{align*}
     \de(\beta_{3\ell+1})&= (101(001)^{\lfloor\frac{\ell-1}{2}\rfloor}(010)^{\lceil\frac{\ell-1}{2}\rceil} 0)^\f \succ (101(001)^{\lfloor\frac{\ell}{2}\rfloor}(010)^{\lceil\frac{\ell}{2}\rceil} 0)^\f=\de(\beta_{3\ell+4}),\\
   \de(\beta_{3\ell+2})&=(101(001)^{\ell-1}00)^\f\succ(101(001)^{\ell}00)^\f=\de(\beta_{3\ell+5}).\end{align*}
  Furthermore,
  \[
  \de(\beta_{3\ell+1})\to 101(001)^\f=\de(\beta_a),\quad \de(\beta_{3\ell+2})\to 101(001)^\f=\de(\beta_a)\quad\textrm{as}\quad \ell\to\f.
  \]
  Therefore, by Lemma \ref{lem:delta-beta} we obtain that
 $\beta_{3\ell+1}\searrow \beta_a$ and $\beta_{3\ell+2}\searrow \beta_a$ as $\ell\to\f$.
\end{proof}

\section*{Acknowledgements}
The first author was supported by NSFC No.~11971079.


\end{document}